\documentclass[12pt]{amsart}

\usepackage{amsmath, amssymb, amscd}
\usepackage[mathscr]{eucal}
\usepackage{mathrsfs}
\usepackage{verbatim}
\usepackage{version}
\usepackage{nicefrac}
\usepackage{pstricks}
\usepackage{pst-all}
\usepackage{color}
\usepackage[all]{xy}
\usepackage{mathdots}

\usepackage{pict2e}

\usepackage[enableskew]{youngtab}
\usepackage{young}

\usepackage{color}
\usepackage{amsmath, latexsym, amssymb}
\usepackage{graphicx}
\usepackage{amsfonts}
\usepackage{subfigure}
\usepackage{amsthm}
\usepackage{amsmath}
\usepackage{graphics}
\usepackage{amssymb}
\usepackage{latexsym}
\usepackage{amscd}
\usepackage{color}
\usepackage{graphicx}
\usepackage{pict2e}

\usepackage[colorlinks=true,citecolor=red,linkcolor=blue]{hyper ref}
\usepackage{lineno}

\newlength\cellsize 
\savebox2{%
\begin{picture}(15,15)
\put(0,0){\line(1,0){15}}
\put(0,0){\line(0,1){15}}
\put(15,0){\line(0,1){15}}
\put(0,15){\line(1,0){15}}
\end{picture}}
\newcommand\cellify[1]{\def\thearg{#1}\def\nothing{}%
\ifx\thearg\nothing
\vrule width0pt height\cellsize depth0pt\else
\hbox to 0pt{\usebox2\hss}\fi%
\vbox to 15\unitlength{
\vss
\hbox to 15\unitlength{\hss$#1$\hss}
\vss}}
\newcommand\tableau[1]{\vtop{\let\\=\cr
\setlength\baselineskip{-16000pt}
\setlength\lineskiplimit{16000pt}
\setlength\lineskip{0pt}
\halign{&\cellify{##}\cr#1\crcr}}}
\savebox3{%
\begin{picture}(15,15)
\put(0,0){\line(1,0){15}}
\put(0,0){\line(0,1){15}}
\put(15,0){\line(0,1){15}}
\put(0,15){\line(1,0){15}}
\end{picture}}
\newcommand\expath[1]{%
\hbox to 0pt{\usebox3\hss}%
\vbox to 15\unitlength{
\vss
\hbox to 15\unitlength{\hss$#1$\hss}
\vss}}

\newcommand{\bM}{\mathbf M}
\newcommand{\bN}{\mathbf N}
\newcommand{\bO}{\mathbf O}

\newcommand{\bV}{\mathbf V}

\newcommand{\cF}{\mathcal F}

\newcommand{\cM}{\mathcal M}
\newcommand{\cN}{\mathcal N}

\newcommand{\cO}{\mathcal O}

\newcommand{\cT}{\mathcal T}
\newcommand{\cU}{\mathcal{U}}
\newcommand{\cV}{\mathcal V}

\newcommand{\fg}{\mathfrak g}
\newcommand{\fh}{\mathfrak h}

\newcommand{\fo}{\mathfrak o}
\newcommand{\fp}{\mathfrak p}

\newcommand{\heta}{{\widehat \eta}}

\newcommand{\of}{\overline{ f}}

\newcommand{\ocO}{\overline{\cO}}
\newcommand{\opi}{\overline{\pi}}

\newcommand{\rO}{\mathrm{O}}

\newcommand{\sA}{\mathsf A}

\newcommand{\sE}{\mathsf E}

\newcommand{\sT}{\mathsf T}

\newcommand{\tmu}{\tilde{\mu}}

\newcommand{\tcT}{\widetilde{\cT}}

\newcommand{\aaa}{\mathbb A}

\newcommand{\pp}{\mathbb P}

\newcommand{\cc}{\mathbb C}
\newcommand{\rr}{\mathbb R}
\newcommand{\zz}{\mathbb Z}

\newcommand{\Coeff}{\mathrm{Coeff} }
\newcommand{\Coker}{\mathrm{Coker} }

\newcommand{\dual}{{^\vee} }

\newcommand{\End}{\mathrm{End} }

\newcommand{\Ext}{\mathrm{Ext} }

\newcommand{\gl}{\mathfrak{gl} }

\newcommand{\GL}{\mathrm{GL} }

\newcommand{\git}{/\!\!/ }
\newcommand{\Hom}{\mathrm{Hom} }

\newcommand{\Ker}{\mathrm{Ker} }
\newcommand{\Image}{\mathrm{Im} }
\newcommand{\iso}{\mathrm{iso} }
\newcommand{\Gr}{\mathrm{Gr} }
\newcommand{\Lie}{\mathrm{Lie} }
\newcommand{\inv}{^{-1}}
\newcommand{\Id}{\mathrm{Id} }
\newcommand{\Map}{\mathrm{Map}}

\newcommand{\mindeg}{\mathrm{min.deg}\,}
\newcommand{\modality}{\mathrm{mod}}
\newcommand{\nilp}{\mathrm{nilp}}

\newcommand{\rank}{\mathrm{rank} }
\newcommand{\reg}{\mathrm{reg} }

\newcommand{\Res}{\operatornamewithlimits{Res}}

\newcommand{\Spec}{\mathrm{Spec}\, }

\newcommand{\fsl}{\mathfrak{sl} }

\newcommand{\sm}{\mathrm{sm}}
\newcommand{\SL}{\mathrm{SL}}
\newcommand{\SO}{\mathrm{SO}}
\newcommand{\Sp}{\mathrm{Sp}}

\newcommand{\Spin}{\mathrm{Spin}}
\newcommand{\SU}{\mathrm{SU}}
\newcommand{\Supp}{\mathrm{Supp}}

\newcommand{\fsp}{\mathfrak{sp} }

\newcommand{\tr}{\mathrm{tr}}

\newcommand{\USp}{\mathrm{USp} }

\newtheorem{prop}{Proposition}[section]
\newtheorem{thm}[prop]{Theorem}
\newtheorem{lem}[prop]{Lemma}
\newtheorem{cor}[prop]{Corollary}
 
\theoremstyle{remark}
 \newtheorem{rk}[prop]{Remark}
  
\theoremstyle{definition}

\numberwithin{equation}{section}
\numberwithin{figure}{section}

\begin{document}


\title{Geometry of Uhlenbeck partial compactification of orthogonal instanton spaces and the K-theoretic Nekrasov partition functions}

\author{Jaeyoo Choy}

\address{Dept. Math., Kyungpook Nat'l Univ., Sangyuk-dong, Buk-gu, Daegu 702-701, Korea}

\email{choy@knu.ac.kr, jaeyoochoy@gmail.com}

\dedicatory{On the occasion of Professor Hiraku Nakajima's fifty-fifth birthday}

\thanks{The author was partially supported by JSPS Ronpaku fellowship.}


\begin{abstract}
Let $\cM^K_n$ be the moduli space of framed $K$-instantons over $S^4$ with instanton number $n$ when $K$ is a compact simple Lie group of classical type.
Let $\cU^{K}_{n}$ be the Uhlenbeck partial compactification of $\cM^{K}_{n}$.
A scheme structure on $\cU^{K}_{n}$ is endowed by Donaldson as an algebro-geometric Hamiltonian reduction of ADHM data.
In this paper, for $K=\SO(N,\rr)$, $N\ge5$,  we prove that $\cU^{K}_{n}$ is an irreducible normal variety with smooth locus  $\cM^{K}_{n}$.
Hence, together with the author's previous results on $\USp(N)$, the K-theoretic Nekrasov partition function for any simple classical group other than $\SO(3,\rr)$, is interpreted as a generating function of Hilbert series of the instanton moduli spaces.

Using this approach we also study the case $K=\SO(4,\rr)$ which is the unique semisimple but non-simple classical group.
\end{abstract}

\maketitle
\setcounter{tocdepth}{1}
\tableofcontents

\thispagestyle{empty} \markboth{Jaeyoo Choy}{Uhlenbeck spaces of orthogonal instantons and the K-theoretic Nekrasov partition functions}


\section{Introduction}
\label{sec: intro}

Let $K$ be a compact classical group.
The K-theoretic Nekrasov partition function is defined by Nekrasov and Shadchin as the formal sum of the equivariant integrations of K-theory classes on the ADHM quiver representation space associated to $K$-instantons \cite{NS}.
Our previous result \cite{Choy} says that it is the generating function of the Hilbert series of the coordinate rings of the framed $K$-instanton moduli spaces $\cM^{K}_{n}$ over all instanton numbers $n$ if $K=\USp(N/2)$ the real symplectic group where $N\in 2\zz_{\ge0}$.
We aim to prove the parallel result for $K=\SO(N,\rr)$.
Note that such a result has been known for $K=\SU(N)$ essentially due to Crawley-Boevey \cite{CB}.


\subsection{Main result}\label{subsec: intro main result}

We state the main result of the paper in this subsection.
For the purpose we need to describe the ADHM data of instantons following Donaldson's argument \cite{Do}.
First we consider the vector space $\bM$ of ordinary ADHM quiver representations coming from framed $\SU(N)$-instantons with instanton number $k$.
It is given as $\bM=\End(V)^{\oplus2}\oplus \Hom(W,V)\oplus \Hom(V,W)$.
$\bM$ is a cotangent space, hence naturally a symplectic vector space.
The adjoint $\GL(V)$-action preserves the symplectic structure.

We fix an instanton number $n$.
Let $k=2n,4n$ or $n$ according to $K=\SO(N,\rr)$ ($N\ge4$), $\SO(3,\rr)$ or $\USp(N/2)$.
Let
    $$
    G:=\left\{\begin{array}{lll} \Sp(n) && \mbox{if $K=\SO(N,\rr)$ with $N\ge4$}
    \\
     \Sp(2n) && \mbox{if $K=\SO(3,\rr)$}
     \\
    \rO(n) && \mbox{if $K=\USp(N/2)$ with $N\in 2\zz_{\ge0}$.}
    \\
    \end{array}
    \right.
    $$
Let $V,W$ be the vector representations of $G,K_{\cc}$ respectively.

We define a symplectic subspace of $\bM$ as
    $$
    \bN:=\bN_{V,W}:=\{(B_1,B_2,i,j)\in \bM|\, B_1=B_1^*,B_2=B_2^*, j=i^*\}
    $$
where $B_1^*,B_2^*,i^*$ denote the right adjoints of $B_{1},B_{2},i$ respectively (see \S\ref{subsec: adj}).
Since the adjoint $G$-aciton preserves the symplectic structure of $\bN$, there is a moment map $\mu$ (see \S\ref{subsec: moment maps}):
        $$
        \mu\colon \bN\to \Lie(G),\quad (B_1,B_2,i,j)\mapsto [B_1,B_2]+ij.
        $$
We call the quiver representations in $\mu\inv(0)$ as SO-\textit{data} or Sp-\textit{data} if $K=\SO(N,\rr)$ or $\USp(N/2)$ respectively.

An element $x=(B_1,B_2,i,j)$ of $\bM$ is defined to be \textit{stable} (resp.\ \textit{costable}) if there is no proper $B_1,B_2$-invariant subspace in $V$ containing $i(W)$ (resp.\ there is no nonzero $B_1,B_2$-invariant subspace contained in $\Ker(j)$).
The \textit{regular} locus $\bM^\reg$ means the locus of stable-costable quiver representations in $\bM$.
For any subset $A\subset \bM$, $A^{\reg}$ denotes $A\cap \bM^{\reg}$.

\begin{thm}\label{th: main}
For $K=\SO(N,\rr)$, $N\ge4$, the moment map $\mu$ is flat.

Moreover the following assertions are true:

$(1)$ Let $N=4$.
Then $\mu\inv(0)$ is a reduced variety with the $n+1$ irreducible components.
And the regular locus $\mu\inv(0)^\reg$ is Zariski dense open in $\mu\inv(0)$.

$(2)$ Let $N\ge5$.
Then $\mu\inv(0)$ is an irreducible normal variety.
And the regular locus $\mu\inv(0)^\reg$ is Zariski dense open in $\mu\inv(0)$.
\end{thm}

In the next subsection we deduce from Theorem \ref{th: main} a scheme-theoretic description of the instanton moduli space $\cM^{K}_{n}$ and the Uhlenbeck space $\cU^{K}_{n}$.


\subsection{Geometry of the instanton moduli space}\label{subsec: Geometry of the instanton moduli space}

We give scheme structures on $\cM^{K}_{n}$ and $\cU^{K}_{n}$ first and then discuss the geometric structures induced from Theorem \ref{th: main}.
Let $K=\SO(N,\rr)$, $N\ge 5$.
Note that $K$ is a simple group.
A theorem of Donaldson \cite[p.459]{Do} says that $\cM^{K}_{n}$ is identified with the quasi-affine quotient $\mu\inv(0)^{\reg}/G$.
Secondly the Uhlenbeck space $\cU^{K}_{n}$ is identified with the affine GIT quotient $\mu\inv(0)\git G$.
This identification follows from comparison between the Uhlenbeck stratification of $\cU^{K}_{n}$ and the stratification of $\mu\inv(0)\git G$ given as 
    \begin{equation}\label{eq: stratification of Uhlenbeck space}
    \mu\inv(0)\git G=\bigsqcup_{0\le k'\le k} \mu_{k'}\inv(0)^\reg/\Sp(V_{k'})\times S^{\frac{k-k'}2}\aaa^2.
    \end{equation}
Here $V_{k'}$ is a symplectic vector space of dimension $k'$, $\mu_{k'}$ denotes the moment map on   $\bN_{V_{k'},W}$ with respect to $\Sp(V_{k'})$ and $S^{n'}\aaa^{2}$ denotes the ${n'}{}^{\mathrm{th}}$ symmetric product.
See \cite[Theorem 2.6 (1)]{Choy} for detail.
Therefore $\cM^{K}_{n}$ and $\cU^{K}_{n}$ are realized as quasi-affine and affine schemes respectively.
Due to the scheme structure of $\cU^{K}_{n}$, we call it the scheme-theoretic Uhlenbeck partial compactification (Uhlenbeck space, for short).
The following assertion is an immediate corollary of Theorem \ref{th: main} (2).

\begin{cor}\label{cor: main}
For $\SO(N,\rr)$, $N\ge5$, the Uhlenbeck space $\cU^{K}_n$ is an irreducible normal variety of dimension $2n(N-2)$ with the smooth locus $\cM^{K}_{n}$.
\end{cor}

Let $K=\SO(4,\rr)$.
It is a rather special classical group in the sense that it is a unique non-simple classical group.
The universal covering of $K$ is $\Spin(4)=\USp(1)\times\USp(1)$ and thus the topological type of a $K$-instanton is given by a pair $(n_{1},n_{2})$ of the instanton numbers of two $\USp(1)$-instantons.
Such a pair $(n_{1},n_{2})$ is also called instanton number (of an $\SO(4,\rr)$-instanton).
We denote by $\cM^{K}_{(n_{1},n_{2})}$ the $K$-instanton moduli space for an instanton number $(n_{1},n_{2})$ of instanton numbers.
It is naturally identified with the product $\cM^{\USp(1)}_{n_{1}}\times\cM^{\USp(1)}_{n_{2}}$.
A partial compactification of $\cM^{K}_{(n_{1},n_{2})}$ can be given as $\cU^{\USp(1)}_{n_{1}}\times \cU^{\USp(1)}_{n_{2}}$. 
Therefore $\cM^{K}_{(n_{1},n_{2})}$ and its partial compactification admit respectively the scheme structures from the Uhlenbeck spaces.

Let $n:= n_{1}+n_{2}$.
Since the associated vector bundle of a $K$-instanton with instanton number $(n_{1},n_{2})$ is an orthogonal vector bundle with $c_{2}=2n$, the disjoint union $\bigsqcup_{n=n_1+n_2}\cM^K_{(n_1,n_2)}$ is realized as a quasi-affine scheme $\mu\inv(0)^{\reg}/G$.
Hence we have two scheme structures on $\cM^{K}_{(n_{1},n_{2})}$.
These two scheme structures are isomorphic via the tensor product morphism between the associated framed vector bundles $(\cF_{1},\cF_{2})\mapsto \cF_{1}\otimes\cF_{2}$.
Here $\cF_{l}$, $l=1,2$, is a framed symplectic vector bundle with rank $2$ and the second Chern class $c_{2}=n_{l}$ and thus $\cF_{1}\otimes\cF_{2}$ is a framed orthogonal vector bundle with rank $4$ and $c_{2}=2n$.
However the natural partial compactification $\mu\inv(0)\git G$ differs from the disjoint union of $\cU^{\USp(1)}_{n_{1}}\times \cU^{\USp(2)}_{n_{2}}$ over the instanton numbers $(n_{1},n_{2})$ with $n_{1}+n_{2}=n$.
The following description on the scheme structure of $\mu\inv(0)\git G$ is an immediate corollary of 
Theorem \ref{th: main} (1).

\begin{cor}\label{cor: main 1}
$\mu\inv(0)\git G$ is a reduced variety with the $n+1$ irreducible components.
These irreducible components are the Zariski closures of $\cM^K_{(n_1,n_2)}$ where $n=n_{1}+n_{2}$ with $n_{1}=0,1,...,n$.
\end{cor}

By the stratification \eqref{eq: stratification of Uhlenbeck space}	
	$$
	\mu\inv(0)\git G= \bigsqcup_{n_{1}'+n_{2}'\le n}\cM^{K}_{(n_{1}',n_{2}')}\times S^{n-n_1'-n_2'}\aaa^2,
	$$ 
we have also a set-theoretic description of the irreducible components and their intersection.	
The Zariski closure of $\cM^{K}_{(n_{1},n_{2})}$ in $\mu\inv(0)\git G$ is stratified by $\cM^{K}_{(n_{1}',n_{2}')}\times S^{n-n_1'-n_2'}\aaa^2$ for $n_{1}'\le n_{1}$, $n_{2}'\le n_{2}$.
The intersection of those Zariski closures consists of the common strata.


\subsection{K-theoretic Nekrasov partition function and Hilbert series of instanton moduli spaces}\label{subsec: K-theoretic Nekrasov partition function}

A motivation of the paper arises from Nekrasov-Shadchin's definition of the K-theoretic Nekrasov partition functions \cite{NS}.
We do not give the precise definition here. 
See the original paper \cite{NS} or the expository part \cite[\S\S1.5--1.6]{Choy}.
An essential observation of the K-theoretic Nekrasov partition is that the normality of the Uhlenbeck space $\cU^{K}_{n}$ assures that this partition function becomes the generating functions of Hilbert series of coordinate rings of the instanton spaces $\cM^{K}_{n}$.
This interpretation was done for $K=\USp(N/2)$ in the author's previous work \cite{Choy}.
We also obtain a similar interpretation for $\SO(N,\rr)$, $N\ge5$ due to the normality result in Theorem \ref{th: main} (2).
A similar argument also applies to the case $K=\SO(N,\rr)$ for $N\ge4$.
But $\mu\inv(0)\git G$ is not normal, so the parallel interpretation fails.
We get back to this point later in this subsection.

There are naturally isomorphic instanton moduli spaces with different structure groups (see \cite[\S1.4]{Choy}).
Such pairs are precisely 
	$$
	\begin{aligned}
	&
	(\SU(2),\USp(1)),\ (\SU(2),\SO(3,\rr)),\ (\USp(1)\times \USp(1),\SO(4,\rr)),
	\\
	&
	(\USp(2),\SO(5,\rr)),\ (\SU(4),\SO(6,\rr)).
	\end{aligned}
	$$ 
For each pair, ADHM descriptions of instantons are different.
Here the ADHM data for $\USp(1)\times\USp(1)$ mean pairs of the ADHM data for $\USp(1)$.
The K-theoretic Nekrasov partition functions do \textit{not} depend on the ADHM descriptions except the pairs $(\SU(2),\SO(3,\rr))$, $(\USp(1)\times\USp(1),\SO(4,\rr))$.
This fact follows from that the scheme structures of the instanton spaces for the pairs $(\SU(2),\USp(1)), (\USp(2),\SO(5,\rr)),(\SU(4),\SO(5,\rr))$ are naturally mutually isomorphic, and that the above interpretation that the K-theoretic partition function is the generating function of the coordinate rings of the instanton spaces. 
See \cite[\S1.4]{Choy} for these natural isomorphisms.

The cases of pairs $(\SU(2),\SO(3,\rr))$, $(\USp(1)\times\USp(1),\SO(4,\rr))$ are rather complicated.
For the former pair $(\SU(2),\SO(3,\rr))$, the affine GIT quotient $\mu\inv(0)\git G$ for $\SO(3,\rr)$ is not even irreducible in general.
In contrast, $\cU^{\SU(2)}_{n}$ is an irreducible normal variety \cite{CB}.
See \cite[Theorem 4.3]{Choy}.

For the latter pair $(\USp(1)\times\USp(1),\SO(4,\rr))$, $\mu\inv(0)\git G$ for $\SO(4,\rr)$ is connected so the Hilbert series has the constant term $1$.
On the other hand the Zariski closure of $\cM^{K}_{(n_{1},n_{2})}$ in the ADHM space for $\USp(1)\times\USp(1)$ becomes the product of Uhlenbeck spaces $\cU^{\USp(1)}_{n_{1}}\times\cU^{\USp(1)}_{n_{2}}$.
Hence we have many connected components according to the decomposition $n=n_{1}+n_{2}$.
For each pair $(n_{1},n_{2})$ with $n=n_{1}+n_{2}$, the Hilbert series of $\cM^{K}_{(n_{1},n_{2})}$ has the constant term $1$.
Thus the Hilbert series for the pairs do not coincide.
In other words the K-theoretic partition function for $\SO(4,\rr)$ does not coincide with the square of the partition function for $\USp(1)$.
At this moment we do not know the precise geometric meaning of the difference between the two partition functions in both cases.


\subsection{Further motivation}

The K-theoretic Nekrasov partition function is a counterpart in the Grothendieck groups of coherent sheaves, of the Nekrasov partition function which is originally defined topologically \cite{Nek}: equivariant integration of the unit class in the in the equivariant cohomology ring of the ordinary Gieseker  spaces of framed torsion-free sheaves when $K=\SU(N)$.
A purpose of the latter topological partition function in \cite{Nek} is to understand the $4$-dimensional $\cN=2$ supersymmetric gauge theory in physics related to the Seiberg-Witten prepotential \cite{SW}.
The lowest degree coefficient of its logarithm turns out to be the Seiberg-Witten prepotential as was proven by Nakajima-Yoshioka \cite{NY} and Nekrasov-Okounkov \cite{NO} independently.
For other compact groups it can be defined via equivariant integration over the Uhlenbeck spaces given by Braverman-Finkelberg-Gaitsgory \cite{BFG}.
In the case, a similar relation between the Nekrasov partition function and the Seiberg-Witten prepotential is shown by Braverman-Etingof \cite{BE}. 

The K-theoretic Nekrasov partition function is proposed as a non-perturbative solution of $5$-dimensional $\cN=1$ gauge theory by Nekrasov \cite{Nek96}.
This $5$-dimensional gauge theory itself is interesting and appears in various other situations in mathematics.
When $K=\SU(N)$, the $F$-terms give rise to the topological string partition function defined on a singular local Calabi-Yau via geometric engineering in physics terminology.
To make precise mathematically, we consider first the local Calabi-Yau $3$-fold $\cO_{\pp^1}(-1)^{\oplus2}$ fibred over $\pp^{1}$.
By taking the quotient by the finite cyclic $\zz_{N}$-action $l\mapsto (\exp(2\pi l\sqrt{-1}/N),\exp(-2\pi l\sqrt{-1}/N))$ on the fibres $\cc^{2}$, we obtain the singular Calabi-Yau $3$-fold $\cO_{\pp^1}(-1)^{\oplus2}/\zz_N$ mentioned above.
Now the above statement in physics becomes the following: 
The K-theoretic Nekrasov partition function $Z^{K}$ after substitution $\hbar=q_1\inv=q_2\inv$ and multiplication by $\hbar^{-2}$, becomes the generating function of Gromov-Witten invariants on the crepant resolution of $\cO_{\pp^1}(-1)^{\oplus2}/\zz_N$.
This was checked by Iqbal and Kashani-Poor  \cite{IKP}. 

There is a perspective from Intriligator-Seiberg's mirror symmetry in $3$-dimensional $\cN=4$ theory in physics \cite{IS}.
It says a certain duality between the two moduli spaces of ALE gravitational instantons and ADE instantons.
In particular the K-theoretic Nekrasov partition function $Z^{K}$ corresponds to the Coulomb branch $\cM_{C}$ of $3$-dimensional $\cN=4$ theory via the Kaluza-Klein reductions (see Benini-Tachikawa-Xie \cite{BTX}).
As an application of Theorem \ref{th: main}, $Z^{K}$ computes the mirror of $\cM_{C}$ for the type $D$ quiver gauge theory.
See a mathematical exposition of this theory due to Nakajima \cite{Nak_Coulomb}\cite{Nak_Questions} and Braverman-Finkelberg-Nakajima \cite{BFN1}\cite{BFN2}.

Let us look at the aspect of the actual computation of $Z^{K}$, which is now understood as the Hilbert series for any classical groups $K$ except $\SO(4,\rr),\SO(3,\rr)$ due to \S\ref{subsec: Geometry of the instanton moduli space}.
For $K=\SU(N)$ there are two ways of computation of $Z^{K}$: (1) a direct computation using Weyl's integral formula \cite[App.\ A]{FL} (cf.\ \cite{EG}), (2) a closed formula a.k.a.\ the blowup equation \cite[Theorem 2.4]{NYII}.
For $K=\SO(N,\rr),\USp(N/2)$, $Z^{K}$ is computed for small instanton numbers \cite{BHM}.
In the formulation of the blowup equation for $K=\SU(N)$ in \cite{NYII}, it is essential to use two Gieseker spaces: one is the moduli space of framed torsion-free sheaves on $\pp^{2}$ while the other is defined over the blown-up $\pp^{2}$.
Nakajima and Yoshioka deduced the blowup equation by comparing the two formulas arising from the localization (with respect to the torus-action on the framings) and the push-forward of some line bundle defined on the Gieseker space over the blown-up $\pp^{2}$. 
For other simple groups $K$, they also gave conjectures \cite[(9.2), (9.3)]{NYI} on the geometry of the Uhlenbeck spaces, under which the blowup equation holds again (cf.\ \cite[(6.10)--(6.14)]{NYI}).
In a similar context, the topological Nekrasov partition function for an arbitrary simple group $K$ turns out to be approximately an eigenfunction of a Schr\"odinger operator due to Braverman-Etingof (\cite[\S4]{BE}, cf.\ \cite[\S3.6]{Braverman}).


\subsection{Flatness of moment map}
\label{subsec: flatness of moment map}

The technical core in this paper is flatness of the moment map $\mu$ defined on $\bN$ for $\SO(N,\rr)$-instantons with rank $N\ge4$.
We remark that $\mu$ is \textit{not} flat if $N=2$ (\cite[Theorem 4.2 (1)]{Choy}).

A first observation is that $\bN$ itself is a cotangent bundle if $N$ is even.
In the case the $G$-action is induced from the base of the cotangent bundle and thus it is Hamiltonian.
The moment map $\mu$ coincides with the canonical one up to a linear isomorphism of $\bN$ and thus the flatness of these two moment maps are equivalent.
Recall that in the cases $K=\SU(N),\USp(N/2)$ the corresponding spaces $\bM,\bN$ are always cotangent spaces and the moment maps are induced by the action on the base.
Thus our explanation from now on equally works for these cases.

In general the canonical moment map $\mu$ on the cotangent bundle $T^*\cV$ is \textit{not} flat, where $\cV$ is a $G$-module.
E.g., the cases $\bM$ and $\bN$ for $K=\SU(0)$ and $\SO(0,\rr),\SO(1,\rr),\SO(2,\rr)$ (see \cite{Knop}\cite{Losev} for more general flatness results).
In fact flatness of $\mu$ amounts to a dimension bound of the zero fibre $\dim \mu\inv(0)\le 2\dim \cV-\dim G$.
By a theorem of Vinberg, if $G$ is a reductive group, $\dim \mu\inv(0)$ is equal to $\dim \cV+\max_{\cV'} \mathrm{tr.deg}(\cc(\cV')^G)$, where $\cV'$ runs over the $G$-invariant locally closed subsets of $\cV$ and $\mathrm{tr.deg}$ stands for transcendence degree (\cite{Vinberg}; see also \cite[\S2.2]{Pa}).
The latter summand $\max_{\cV'} \mathrm{tr.deg}(\cc(\cV')^G)$ is called \textit{modality} of $\cV$ and in our paper its definition is slightly different but equivalent (see \S\ref{subsec: modality}).
So the flatness of $\mu$ amounts to the inequality $\max_{\cV'} \mathrm{tr.deg}(\cc(\cV')^G)\le \dim\cV-\dim G$.

The above modality inequality for flatness for the space of symmetric pairs was used by Panyushev \cite{Pa} (see also Brennan \cite{Br}).
It corresponds to the rank 0 case of USp-instantons.
The flatness for USp-instantons with the higher ranks comes from the base change argument \cite{Choy}.
For $K=\SU(N),\SO(N,\rr)$ the flatness does \textit{not} any more follow from the first factor $\gl(V)$ or $\fp(V)$ solely, but also needs the second factor $\Hom(W,V)$.
In a similar context Crawley-Boevey used the modality inequality to determine the flatness of the moment maps defined on the space of representations of the double of a general quiver \cite[\S3]{CB1}.

We introduce a new aspect and approach used in the proof of flatness for $\SO(N)$-data ($N\ge4$) compared to the above known cases $K=\SU(N),\USp(N/2)$.
Here $\SO(N)$-data mean the ADHM data for $\SO(N,\rr)$-instantons.
In the even rank case, the base $\cV$ of the cotangent space $T^*\cV$ can be written in the form $\cV_1\oplus\cV_2$ where $\cV_{1}\subset \End(V)$ and $\cV_{2}\subset \Hom(W',V)$ for a maximal isotropic subspace $W'$ in $W$ (see \S\ref{subsec: moment maps}).
In fact the flatness of $\mu_{\cV}$ with respect to $G$ amounts to the flatness of $\mu_{\cV_{2}}$ with respect to the stabilizer subgroup $G^{x}$ of a generic nilpotent endomorphism $x\in\cV_1$ due to a reciprocal property of modality (see \S\ref{subsec: reduction lemma}).
So the original flatness problem is replaced to the modality of $\cV_{2}$ with respect to $G^{x}$.
For the known cases $K=\SU(N),\USp(N/2)$, this modality is rather easy to compute since $G^{x}$ is abelian (see \S\ref{subsec: reduction lemma for ADHM data}).
In contrast, if $K=\SO(N,\rr)$ with $N\in2\zz_{\ge2}$, $G^{x}$ is non-abelian.
For instance if $x$ is a generic nilpotent element, $\Lie(G^{x})$ is isomorphic to the truncated current algebra $\fsl_{2}[z]/z^{n}\fsl_{2}[z]$ (Lemma \ref{lem: fg}).
Moreover $\cV_{2}=\cc^{N/2}\otimes \cc^{2}\otimes \cc[z]/(z^{n})$, where $\cc^{2}\otimes \cc[z]/(z^{n})$ is the vector representation of $\fsl_{2}[z]/z^{n}\fsl_{2}[z]$.
We compute the modality of $(G^{x},\cV_{2})$ (Theorem \ref{th: flat half loop algebra}), which assures the flatness of the moment map.

Note that in the case of $\SO(3)$-data, the above method does \textit{not} work as $\bN$ is not a cotangent space.
Hence we will use a different approach to this problem in \cite{Choy2}.


\subsection{Contents}
\label{subsec: intro contents}

This paper is organized as follows.
In \S\S\ref{subsec: symmetric and antisymmetric}--\ref{subsec: factorization}, we first review the basic materials on moment maps $\mu$ for $\SO(N)$-data.
Then in \S\ref{subsec: the proof of main theorem}, we prove Theorem \ref{th: main} assuming the flatness of $\mu$ for $N\ge4$ (Lemma \ref{lem: reduction}) and the normality of $\mu\inv(0)$ in the instanton number $n=2$ case and $N\ge5$ (Lemma \ref{lem: normal k4N5}).
In \S\ref{subsec: further description on irreducible components}, we give a quiver description of each  irreducible components of $\mu\inv(0)$ when $N=4$.
This will be done by giving a quiver-theoretic description of the tensor product morphism of framed vector bundles (Theorem \ref{th: ADHM tensor product morphism}).

In \S\ref{sec: flatness via modality}, we prove the flatness of $\mu$.  
In \S\S\ref{subsec: modality}--\ref{subsec: flatness}, we define modality of a Hamiltonian vector space and then deduce its basic properties, e.g.\ an inequality of modality equivalent to the flatness of the moment map.
In \S\ref{subsec: r2k2}, using the modality inequality we show the flatness of the moment map on a vector representation with respect to  a truncated current algebra.
Using this flatness, we show the flatness of the moment map for $\SO(N)$-data in \S\ref{subsec: reduction lemma}.

In \S\ref{sec: normality}, we prove normality of $\mu\inv(0)$ for $N\ge5$.
By the base change argument and the factorization property, our study reduces to the case $N=5,n=2$. 
The method is based on Kraft-Procesi's theory on nilpotent pairs.

In \S\ref{sec: reduction lemma for ADHM}, we apply the argument similar to the known cases for the flatness and the normality: the ordinary ADHM data and the Sp-data. 
As a by-product we obtain a generalization of a result of Gan-Ginzburg \cite[Proposition 2.3.2]{GG}.

\vskip.3cm \noindent\textit{Acknowledgement.}
The author would like to express his deepest gratitude to the supervisor Professor Hiraku Nakajima who has guided this topic and related subjects.
The argument in this paper has been improved thanks to his suggestion on the factorization property (\cite{BFN}\cite{MO}) from earlier versions. 

Most part of this research was done during the author's JSPS Ronpaku fellowship 2010--2014 at Research Institute for Mathematical Sciences, Kyoto University.
The author is grateful to hospitality from the staffs there.
He is also grateful to Professors Hoil Kim, Bumsig Kim for their encouragement and hospitality at Korea Institute for Advanced Study where a part of the paper is written.
He appreciates Dr.\ Ada Boralevi's interest and comments on this work.


\section{Moment maps $\mu$ for SO-data and the proof of Theorem \ref{th: main}}
\label{sec: preliminary}

In this section we reformulate the moment maps $\mu$ for $\SO(N)$-data, $N\in 2\zz_{\ge2}$, so that it is defined over a cotangent space, and then prove the main theorem (Theorem \ref{th: main}).
We also describe the tensor product isomorphism between an irreducible component of $\mu\inv(0)^{\reg}/G$ and $\cM^{\SO(4,\rr)}_{(n_{1},n_{2})}$ in terms of quiver representations.

In \S\S\ref{subsec: symmetric and antisymmetric}--\ref{subsec: trace pairing}, we give basic formulation and properties of the space of $\SO$-data.
In \S\ref{subsec: moment maps}, we give the aforementioned reformulation of the moment map $\mu$.
In \S\S\ref{subsec: reduction}--\ref{subsec: factorization}, as preliminary steps toward the proof of Theorem \ref{th: main}, we study a stronger flatness assertion for $\mu$ (Lemma \ref{lem: reduction}) and a factorization property of the zero scheme $\mu\inv(0)$ (Lemma \ref{lem: quotient}).
In \S\ref{subsec: the proof of main theorem}, we give the proof of the theorem.
In \S\ref{subsec: further description on irreducible components}, we identify the irreducible components of $\mu\inv(0)$ via the tensor product morphism.
For the purpose we use commutativity of the tensor product morphism and the factorization morphism. 
By this commutativity we will obtain a quiver description of the tensor product morphism (Theorem \ref{th: ADHM tensor product morphism}).

The conventions in the entire part of this paper are as follows.
We are working over $\cc$.
Schemes are of finite type and vector spaces are finite dimensional.
Morphisms are called irreducible, normal, equi-dimensional, Cohen-Macaulay if all the nonempty fibres are so, respectively.
Dimension of a scheme of finite type is defined to be the maximum of dimensions of irreducible components.


\subsection{Symmetric and antisymmetric forms}
\label{subsec: symmetric and antisymmetric}

Let $V$ be a vector space with a nondegenerate form $(\,,\,)_V$.
We denote by $G(V)$ the subgroup of $\GL(V)$ preserving $(\,,\,)_V$.
Let $\fg(V):=\Lie(G(V))$.


We define the space of symmetric forms:
     $$
     \fp(V):=\{A\in \gl(V)|\, (Av,w)_V=(v,Aw)_V\ \mbox{for all $v,w\in V$}\},
     $$
Note that $\fg(V)$ is the space of antisymmetric forms.

In this section we drop the notation $V$ from $\fp(V),\fg(V),\gl(V),\fsl(V)$, etc.\ if no confusion arises.

For the commutator bracket $[A,B]:=AB-BA$, we have the following relations by direct calculation:
    \begin{equation}\label{eq: ftp}
    [\fg,\fp]\subset \fp,\ \ [\fp,\fp]\subset \fg.
    \end{equation}
For the anti-commutator bracket $[A,B]_+:=AB+BA$, we have the following relations by direct calculation:
    \begin{equation}\label{eq: ftp+}
    [\fg,\fg]_+\subset \fp, \ \ [\fg,\fp]_+\subset \fg,\ \ [\fp,\fp]_+\subset \fp.
    \end{equation}


\subsection{Adjoint}\label{subsec: adj}

Let $V_1,V_2$ be vector spaces.
Let $L(V_1,V_2):=\Hom(V_1,V_2)$ and $E(V_1):=\End(V_1)$.

Suppose that $V_1,V_2$ have nondegenerate bilinear forms $(\,,\,)_{V_1},(\,,\,)_{V_2}$ respectively.
For $i\in L(V_1,V_2)$, we define the right adjoint $i^*\in L(V_2,V_1)$, i.e.,  $(v,i^*w)_{V_1}=(iv,w)_{V_2}$.

The map $*\colon L(V_1,V_2)\to L(V_2,V_1),\   i\mapsto i^*$ is a $\cc$-linear isomorphism.
Further if $V_3$ is another vector space with a nondegenerate bilinear form, we have $(ji)^*=i^*j^*$ where $i\in L(V_1,V_2),\ j\in L(V_2,V_3)$.

Let us consider the composite $**$ which is an endomorphism of $L(V_1,V_2)$.
We have
    $$
    **=\left\{\begin{array}{lll}
        \Id& & \mbox{if both $V_1,V_2$ are either symplectic or orthogonal}
        \\
        -\Id & & \mbox{if one of $V_1,V_2$ is symplectic and the other is orthogonal.}
        \end{array}
        \right.
    $$
Hence in the case $**=-\Id$, we have $i^*i\in \fg(V_1)$ and $ii^*\in \fg(V_2)$. Similarly in the case $**=\Id$, we have $i^*i\in \fp(V_1)$ and $ii^*\in \fp(V_2)$.


\subsection{Trace pairing}\label{subsec: trace pairing}

Suppose that $V$ is either symplectic or orthogonal.

Let us consider the trace pairing $\tr\colon \fg\times \fg\to \cc$ the restriction of the ordinary trace pairing for $\gl$.
We claim it is nondegenerate.
As the trace pairing for $\gl$ is nondegenerate, for any $A\in \fg$ there is $B\in \gl$ such that $\tr(AB)\neq0$.
We have
    $$
    \tr(AB)= \tr\left(A\left(\frac{B+B^*}2+\frac{B-B^*}2\right)\right).
    $$
We observe that $\frac{B+B^*}2\in \fp\ \mbox{and}\ \frac{B-B^*}2\in\fg$.
So the claim follows from $\tr (AC)=0$ whenever $C\in \fp$.
We notice that $\tr(AC)=\frac12\tr([A,C]_+)$ and $[A,C]_+\in \fg$ by \eqref{eq: ftp+}.
Thus $\tr(AC)=0$ as  $\fg\subset \fsl$.

The trace pairing $\tr\colon \fp\times \fp\to \cc$ is also nondegenerate by a similar proof.


\subsection{Moment maps}
\label{subsec: moment maps}

A purpose of the subsection is to identify the moment map for the SO-data with even rank $\ge4$ with a natural moment map on a cotangent space.

Let us recall the definition of moment maps.
Let $(V,\omega)$ be a symplectic vector space with a Hamiltonian group action by an algebraic group $G$.
Then a map $\mu\colon V\to \fg\dual$ is called a \textit{moment map} if it satisfies a condition on differential
    $$
    d(\mu)_{x}\colon V\to \fg\dual,\quad v\mapsto \omega(\bullet.x,v).
    $$
It exists uniquely if we impose a constraint $\mu(0)=0$.
It is given by $\mu(x)=\frac12\omega(\bullet x,x)$.

Let $\cV$ be a (linear) $G$-representation.
Then the moment map on the cotangent space $\cV\oplus \cV\dual$ equipped with the standard symplectic form with respect to the diagonal $G$-action is given by  $(x,f)\mapsto f(\bullet.x)$.
We denote this moment map by $\mu_\cV$ to emphasize $\cV$.

However not all Hamiltonian actions on a symplectic vector space arise in this way, e.g.\ the $\SL_2$-action on $\cc^2$.
The same problem lies in the moment map on the space of ADHM quiver representations for $\SO(N)$-data of odd rank $N$.
The moment map cannot be written in the form $\mu_{\cV}$.

Now we identify the moment maps defining SO-data with $\mu_\cV$ for some $\cV$.
Our argument also works for Sp-data with any rank.

We consider the space  of usual ADHM quiver representations
        $$
        \bM_{V,W}=E(V)^{\oplus2}\oplus L(W,V)\oplus L(V,W)
        $$
for some vector spaces $V,W$ as in \S\ref{subsec: intro main result}.
Let $\cV:=E(V)\oplus L(W,V)$ and $G:=\GL(V)$.
Then $\bM_{V,W}$ is identified with the cotangent space $T^{*}\cV$ via trace pairing.
The moment map $\mu$ of $\bM_{V,W}$ with respect to $G$ is given by $(B_1,B_2,i,j)\mapsto [B_1,B_2]+ij$.
Here we identified $\fg$ with $\fg\dual$ via trace pairing.
Thus we have $\mu=\mu_{\cV}$.
The elements in the zero locus $\mu\inv(0)$ are called ADHM data.

The space of ADHM quiver representations for SO-instantons (resp.\ Sp-instantons)
        $$
        \bN_{V,W}=
        \left\{ (B_1,B_2,i,j)\in\bM\middle|\, B_1,B_2\in\fp(V), j=i^*\right\}
        $$
is a symplectic subspace of $\bM_{V,W}$
where $(V,W)$ is a pair of symplectic and orthogonal vector spaces (resp.\ orthogonal and symplectic vector spaces).
Here we used that $L(W,V)$ is a symplectic subspace of $T^{*}L(W,V)=L(W,V)\oplus L(V,W)$ via $L(W,V)\to T^{*}L(W,V),\ i\mapsto(i,i^{*})$.
The subgroup $G(V)$ of $G$ acts on $\bN_{V,W}$ and preserves the symplectic form.
Thus the moment map $\mu$ on $\bN_{V,W}$ is the restriction of the moment map on $\bM_{V,W}$ composed with the projection $\Lie(G)\dual\to \Lie(G(V))\dual$.
Recall that the elements in the zero locus $\mu\inv(0)$ have been called SO- or Sp-data.

We use the notations $\bM,\bN$ instead of $\bM_{V,W},\bN_{V,W}$ respectively if there is no confusion.

Now we rewrite $\mu$ defining SO- and Sp-data in terms of $\mu_\cV$ for some $\cV$ when a decomposition $W=W_1\oplus W_2$ is given.
We decompose $L(W,V)=L(W_1,V)\oplus L(W_2,V)$ and $L(V,W)=L(V,W_1)\oplus L(V,W_2)$.
We consider the two cases:
	\begin{itemize}
	\item[(i)] $W_1,W_2$ are maximal isotropic subspaces in $W$ and $(\,,\,)_W\colon W_1\times W_2\to \cc$ is nondegenerate.
	\item[(ii)] $W=W_1\oplus W_2$ is an orthogonal decomposition.
	\end{itemize}

In the case (i), the right adjoint $*$ sends the summand $L(W_1,V)$ (resp.\ $L(W_2,V)$) isomorphically to the summand $L(V,W_2)$ (resp.\ $L(V,W_1)$).
Therefore any element $i=(i_1,i_2)\in L(W,V)=L(W_1,V)\oplus L(W_2,V)$ maps via $*$ to $(i_2^*,i_1^*)\in L(V,W_1)\oplus L(V,W_2)$.
Therefore we have the isomorphism $L(W,V)\stackrel\cong\to L(W_1,V)\oplus L(V,W_1), i\mapsto (i_1,i_2^*)$.
This isomorphism pulls back the standard symplectic form on $T^{*}L(W_1,V)=L(W_1,V)\oplus L(V,W_1)$ to the half of the one on the symplectic subspace $L(W,V)$ of $T^{*}L(W,V)$.

As a consequence, the moment map $\mu$ defining the SO-data on $\bN$ is identified with $\mu_{\cV}$ composed with a linear isomorphism between $\bN$ and $T^{*}\cV$ where $\cV=\fp(V)\oplus L(W_1,V)$.
This means that study of $\mu$ defined on $\bN$ can be replaced by the study of the latter moment map $\mu_{\cV}$ on $T^{*}\cV$.

In the case (ii), the right adjoint $*$ maps $L(W_l,V)$ onto its dual vector space $L(V,W_l)$ for each $l=1,2$.
Hence $\mu$ is the sum of the two moment maps on $\bN_{V,W_1}$ and on $L(W_2,V)$.
We denote these two moment maps by $\mu_1,\mu_2$ respectively.

Here is the strategy for study on the geometry of $\mu\inv(0)$ in the SO-case.
If $W$ has even dimension, we use a decomposition of $W$ by two complementary maximal isotropic subspaces.
Then $\mu$ is identified with $\mu_\cV$ up to linear isomorphism  and thus $\mu\inv(0)$ is isomorphic to $\mu_\cV\inv(0)$.

If $W$ has odd dimension, we use an orthogonal decomposition $W=W_1\oplus W_2$ with $\dim W_2=1$.
Then $\mu$ is the sum $\mu_1+\mu_2$.
Thus $\mu\inv(0)$ is isomorphic to the fibre product of $\mu_1$ and $-\mu_2$ over $\fg\dual$.
Hence the geometry of $\mu\inv(0)$ comes from that of $\mu_1$ by the base change.
E.g., if $\mu_{1}\inv(0)$ is an irreducible normal variety, then so is $\mu\inv(0)$.
See \S\ref{App: normality} for the proof of this assertion.

\begin{rk}\label{rk: simplify moment2}
In the Sp-case, we used the decomposition $W=0\oplus W$ to check that $\mu\inv(0)$ is a normal variety \cite{Choy}.
\end{rk}


\subsection{a Stronger flatness assertion of moment maps for SO-data}\label{subsec: reduction}

The flatness of the moment maps $\mu$ for SO-data comes from a stronger flatness assertion (Lemma \ref{lem: reduction}).
This claim will be proven later (\S\ref{subsec: reduction lemma}), but with this we prove Theorem \ref{th: main} in \S\ref{subsec: the proof of main theorem}.

Let $V,W$ be a symplectic and orthogonal vector spaces of dimensions $k,N$ respectively.
Let
	\begin{equation}\label{eq: bN different identification}
	\bN:=\bN_{V,W}=\fp(V)^{\oplus2}\oplus L(W,V),
	\end{equation}
where we used the identification via the obvious projection.
Let $\mu$ be the moment map defined in $\bN$ for $\SO(N)$-data.

Let $E_{k}\colon \fp(V)\to S^{k}\cc$ be the morphism sending $B$ to the set of eigenvalues counted with multiplicities.
It is nothing but the GIT quotient map by the adjoint action of $\Sp(V)$ onto the image.
Note that since the generalized eigenspaces of any endomorphism $B\in \fp(V)$ are even dimensional (Corollary \ref{cor: G}), the image of $E_{k}$ coincides with $S^{k/2}\cc$ embedded in $S^{k}\cc$ as a closed subscheme by the diagonal map.

\begin{lem}\label{lem: reduction}
Let $N\ge4$.
Then the morphism
	$$
	\tmu \colon\bN\to \fg(V)\times S^{k/2}\cc,\quad x=(B_{1},B_{2},i)\mapsto (\mu(x),E_{k}(B_{1}))
	$$
is flat.
\end{lem}

The proof will appear in \S\ref{subsec: reduction lemma}.
As a corollary of the lemma, $\mu$ is flat for $N\ge4$ by composing the above morphism with the projection to $\fg(V)$.

Let $\fp(V)_{l}:=\{B\in\fp(V)|\mbox{$B$ has $l$ distinct eigenvalues}\}$.
Then $\fp(V)_{\ge l}$ is Zariski open in $\fp(V)$.
Let $\mu\inv(0)_{l}:=\mu\inv(0)\cap p\inv(\fp(V)_{l})$ where $p\colon\bN\to \fp(V)$ is the first projection.

\begin{cor}\label{cor: reduction}
If $1\le l\le k/2$, $\dim\mu\inv(0)_{l}=\dim\bN-\dim\fg(V)-(k/2-l)$.
Otherwise $\mu\inv(0)_{l}$ is empty.
\end{cor}

\proof
Note that $\tmu$ is equi-dimensional by the above lemma and also that $\mu\inv(0)_l=\tmu\inv(E_k(\fp(V)_l))$.
The corollary follows from $\dim E_k(\fp(V)_l)=l$.
\qed\vskip.3cm


\subsection{Local structure of $\mu\inv(0)$ via factorization property}\label{subsec: factorization}

We prove the assertions on $\mu\inv(0)$ in Theorem \ref{th: main} other than the flatness of $\mu$.
We use Drinfeld's factorization property (\cite{BFG}).
This property is formulated for spaces of maps in \cite[Proposition 2.17]{BFG} (see Proposition \ref{prop: Drinfeld's factorization property}).
In this section we follow the quiver-theoretic argument in \cite[Lemma 12.3.2]{MO}.

For any given partition $\eta=(\eta_{1},\eta_{2},...,\eta_{e})$ of $k$, we denote by $S^{\eta}\cc$ the product of $S^{\eta_{l}}\cc$.
Let $(S^{\eta}\cc)_{0}$ be the affine open subset of $S^{\eta}\cc$ consisting of the $e$-tuples $(Z_{1},Z_{2},...,Z_{e})$ such that the supports $\Supp(Z_{l})$ are mutually disjoint.
Let $\Delta S^{\eta}\cc$ be the Zariski closed subset of $S^{\eta}\cc$ consisting of $(Z_{1},Z_{2},...,Z_{e})$ such that $\Supp(Z_{l})$ is a one-point set respectively.
Let $(\Delta S^{\eta}\cc)_{0}:=(S^{\eta}\cc)_{0}\cap \Delta S^{\eta}\cc$.
This is an affine open subset of $\Delta S^{\eta}\cc$.
Let $A_{\eta}\colon S^{\eta}\cc\to S^{k}\cc$, $(Z_{1},Z_{2},...,Z_{e})\mapsto \sum_{l=1}^{e}Z_{l}$.
It is a quasi-finite map and the restriction $A_{\eta}|_{(S^{\eta}\cc)_{0}}$ is an \'etale morphism.
There is a decomposition of $S^{k}\cc$ as
	\begin{equation}\label{eq: decomposition of Skcc}
	S^{k}\cc=\bigsqcup_{|\eta|=k} A_{\eta}((\Delta S^{\eta}\cc)_{0}).
	\end{equation}
Here $|\eta|=\sum_{l=1}^{e}\eta_{l}$.

Let $\pi\colon\bN\to S^{k}\cc$, $(B_{1},B_{2},i)\mapsto E_k(B_{1})$.
Using $\pi,A_{\eta}$, we consider fibre products $\bN\times_{S^{k}\cc}(S^{\eta}\cc)_{0}$ and $ \mu\inv(0)\times_{S^{k}\cc} (S^{\eta}\cc)_{0}$.
We notice that unless all $\eta_{l}$ are even, these fibre products are empty (Corollary \ref{cor: G}).

From now on we may assume $\eta_{l}$ are all even.
We fix mutually orthogonal symplectic subspaces $V_{1}^{0},V_{2}^{0},...,V_{e}^{0}$ in $V$ with $\dim V_l^0=\eta_l$.
Let $\mu_{l}$ be the moment map on $\bN_{V_{l}^{0},W}$ with respect to $\Sp(V_{l}^{0})$ for each $l$.
Let
	\begin{equation}\label{eq: X'Y'}	
 	X':=\left(\bN_{V_{1}^{0},W}\oplus\bN_{V_{2}^{0},W}\oplus\cdots
	\oplus\bN_{V_{e}^{0},W}\right)\times_{S^{\eta}\cc}(S^{\eta}\cc)_{0},
	\quad
	Y':=\bN\times_{S^{k}\cc}(S^{\eta}\cc)_{0}
	\end{equation}
for short.
Note that $X',Y'$ are affine schemes and that $X'$ is Zariski open in $\bN_{V_{1}^{0},W}\oplus\bN_{V_{2}^{0},W}\oplus\cdots\oplus\bN_{V_{e}^{0},W}$.
These facts follow from that $A_{\eta}|_{(S^{\eta}\cc)_{0}}$ is a finite morphism onto the image, hence affine.

Using the orthogonal decomposition $V=V_{1}^{0}\oplus V_{2}^{0}\oplus\cdots\oplus V_{e}^{0}$, we define a morphism
	\begin{equation}\label{eq: pi'eta}
	\textstyle
	\pi_{\eta}'\colon X'\to Y',\quad
 	(B_{l,1},B_{l,2},i_{l})_{1\le l\le e}
	\mapsto (\bigoplus_{l=1}^{e}B_{l,1},B_{2},\bigoplus_{l=1}^{e}i_{l},Z_{1},Z_{2},...,Z_{e}),
	\end{equation}
where the matrix $B_{2}$ is determined so that its block matrix $B_{2}^{(m,l)}\colon V_{l}^{0}\to V_{m}^{0}$ satisfies
 	\begin{equation}\label{eq: factorization}
	\left\{
	\begin{array}{lll}
	B_{2}^{(l,l)}=  B_{l,2} &&
	\\
	B_{m,1}B_{2}^{(m,l)}-B_{2}^{(m,l)}B_{l,1}+i_{m}i_{l}^{*}=0
	&& \mbox{if $m\neq l$}
	\end{array}	
	\right.
	\end{equation}
and $Z_{l}$ is the set of eigenvalues of $B_{1}|_{V_{l}^{0}}$ counted with multiplicities.
Note that $B_{2}^{(m,l)}$ $(m\neq l)$ is determined uniquely by \eqref{eq: factorization} as $B_{m,1}$ and $B_{l,1}$ have different eigenvalues.
Hence each matrix element of $B_2$ is given as a rational function in the matrix elements of $B_{l,1},B_{l,2},i_{l}$ which is defined on $X'$.
Note also that $B_2\in \fp(V)$.
This can seen by taking the right adjoint to each equation in  \eqref{eq: factorization}.
We observe that the second equation of \eqref{eq: factorization} is the non-diagonal block matrix parts of the equation $\mu(\bigoplus_{l=1}^{e}B_{l,1},B_{2},\bigoplus_{l=1}^{e}i_{l})=0$.

Let $\cc[X'],\cc[Y']$ be the affine coordinate rings of the affine schemes $X',Y'$.
Then $\pi_{\eta}'$ induces an algebra homomorphism $(\pi_{\eta}')^{*}\colon \cc[Y']\to \cc[X']$.
We observe that the composite of $\pi_{\eta}'$ and the projection $Y'\to \bN$ is an immersion.
Therefore the following assertion is obvious:

\begin{lem}\label{lem: cl embed}
$(\pi_\eta')^{*}$ is a surjective algebra homomorphism. 
Its kernel is the ideal generated by the matrix element coordinates of $B_{2}^{(l,l)}-B_{l,2}$ and $B_{m,1}B_{2}^{(m,l)}-B_{2}^{(m,l)}B_{l,1}+i_{m}i_{l}^{*}$ where $1\le m\neq l\le k$.

In particular, $\pi_{\eta}'$ is a closed embedding.
\end{lem}

The regular locus $(X')^{\reg}$ maps to $(Y')^{\reg}$ via $\pi_{\eta}'$ due to the following lemma: 

\begin{lem}\label{lem: costable}
Let $x_{l}=(B_{l,1},B_{l,2},i_{l})\in \bN_{V_{l},W}$ for each $l=1,2,...,e$.
Suppose that $x=(B_{1},B_{2},i)\in X'$ where $B_1=\bigoplus_{l=1}^e B_{l,1}$, $B_2^{(l,l)}=B_{l,2}$ and $i=\bigoplus_{l=1}^e i_l$.
If $x_{l}$ is costable for each $l$, $x$ is also costable.

Hence $x$ is regular if so is $x_{l}$ for each $l$.
\end{lem}

\proof
The latter statement is a direct consequence of the former one since costability and stability are equivalent for the SO-data.

We prove the former statement.
Let $K$ be a $B_1,B_2$-invariant subspace of $\Ker(i^*)$.
By $B_1$-invariance we have decomposition $K=\bigoplus_{l=1}^e(K\cap V_l^{0})$.
For, $K$ is decomposed into the generalized eigenspaces of $B_1|_K$.
Each of them is a subspace of $V_l^{0}$ for a precisely one $l$, because $B_{1}|_{V_l^{0}},B_{1}|_{V_{l'}^{0}}$, $l\neq l'$,  do not have a common eigenvalue by the definition of $V_{l}^{0}$.

Therefore $K\cap V_l^{0}$ is a $B_{l,1},B_{l,2}$-invariant subspace of $\Ker(i^*_l)$ for each $l$.
By costability of $x_l$, we have $K\cap V_l^{0}=0$.
\qed\vskip.3cm

\begin{rk}
In the above lemma, we did not use the second constraint of \eqref{eq: factorization} on $x$.
\end{rk}

Let us consider closed subschemes in $X',Y'$ defined as the moment map $0$ loci as:
	$$
	 X:=\left(\mu_{1}\inv(0)\times\mu_{2}\inv(0)\times\cdots\times\mu_{e}\inv(0)\right)\times_{S^{\eta}\cc}(S^{\eta}\cc)_{0},
	\quad
	Y:=\mu\inv(0)\times_{S^{k}\cc}(S^{\eta}\cc)_{0}.
	$$
The composite of $\mu$ with $\pi_{\eta}'$ becomes the restriction of the product morphism $\prod_{l=1}^{e}\mu_{l}$ defined over $X'$.
Thus $\pi'_{\eta}(X)\subset Y$.
We denote by $\pi_{\eta}$ the restriction map from $X$ to $Y$.

We define $G:=\Sp(V)$ and its subgroup $H:=\prod_{l=1}^{e}\Sp(V_{l}^{0})$.
Let $\fg,\fh$ be their Lie algebras respectively.
There are natural $G$-actions on $Y,Y'$ by giving the trivial $G$-action on the factor $S^{\eta}\cc$.
Similarly there are natural $H$-actions on $X,X'$.
We consider another $H$-action on $G\times X'$ by $h.(g,x)=(gh\inv,h.x)$.

\begin{lem}\label{lem: quotient}
$(1)$
$\pi_\eta,\pi_\eta'$ are $H$-equivariant closed embeddings.

$(2)$
The action map $\sigma\colon G\times X\to Y,\ (g,x)\mapsto g.\pi_\eta(x)$ is a surjective smooth morphism with fibres isomorphic to $H$.
\end{lem}

\proof
(1)
We know $\pi_{\eta}'$ is a closed embedding (Lemma \ref{lem: cl embed}).
Recall that $\pi_{\eta}(X)$ is the closed subscheme of $\pi_{\eta}'(X')$ defined by $\mu=0$.
Thus $\pi_{\eta}$ is also a closed embedding.

Let $h=(h_{1},h_{2},...,h_{e})\in H$.
To see the $H$-equivariancy, we check directly 
	$$
	\pi_{\eta}'(h_{l}.(B_{l,1},B_{l,2},i_{l})_{1\le l\le e}) 
	= h.\pi_{\eta}'((B_{l,1},B_{l,2},i_{l})_{1\le l\le e}).
	$$
The only nontrivial part is that  the $(m,l)^{\mathrm{th}}$ block matrix of the second factor, say $C_{2}$, of LHS coincides with $h_{m}B_{2}^{(m,l)}h_{l}\inv$.
By the construction of $\pi_{\eta}'$, $C_{2}$ satisfies 
	$$
	\left\{
	\begin{array}{lll}
	C_{2}^{(l,l)}= h_{l}.B_{l,2} &&
	\\
	(h_{m}.B_{m,1}) C_{2}^{(m,l)}-C_{2}^{(m,l)} (h_{l}.B_{l,1})
	+h_{m}i_{m}i_{l}^{*}h_{l}\inv=0
	&& \mbox{if $m\neq l$.}
	\end{array}	
	\right.
	$$
On the other hand, by multiplying $h_{m}$ from the left and $h_{l}\inv$ from the right to the equations \eqref{eq: factorization}, we get  
	$$
	\left\{
	\begin{array}{lll}
	h_{l}.B_{2}^{(l,l)}= h_{l}.B_{l,2} &&
	\\
	(h_{m}.B_{m,1})(h_{m}B_{2}^{(m,l)}h_{l}\inv)-(h_{m}B_{2}^{(m,l)}h_{l}\inv)(h_{l}.B_{l,1})
	+h_{m}i_{m}i_{l}^{*}h_{l}\inv=0
	&& \mbox{if $m\neq l$.}
	\end{array}	
	\right.
	$$
By comparing these two systems of equations, we obtain $C_{2}^{(m,l)}=h_{m}B_{2}^{(m,l)}h_{l}\inv$ by uniqueness of solution.

(2)
We denote the $G$-action map by $\sigma'\colon G\times X' \to Y'$, $(g,x)\mapsto g.\pi_{\eta}'(x)$.
We claim that $\sigma'$ is a submersion onto the image $\Image(\sigma')$ with the fibres isomorphic to $H$.
Let us finish the proof of the item (2) assuming for a while this claim.
Recall that $\Image(\sigma')$ is a subscheme defined only by the second equation of \eqref{eq: factorization} while $Y$ is defined by the additional equations from $\mu=0$ on the diagonal block matrices.
Therefore $Y$ is a closed subscheme of $\Image(\sigma')$.
It is clear that the pull-back via $\sigma'$ of $\mu=0$ defines $G\times X$. 
Hence $\sigma$ is obtained by the base restriction of the $H$-fibration $\sigma'\colon G\times X'\to \Image(\sigma')$ to $Y$.

We now prove the above claim in the rest of the proof.
First we describe the set-theoretic $\sigma'$-fibres.
For the $H$-action on $G\times X'$ given by $h.(g,x)=(gh\inv,h.x)$, we set $\bO_{(g,x)}$ to be the free $H$-orbit of $(g,x)$ in $G\times X'$.
There is a set-theoretic identification
	$$
	\sigma'{}\inv(\sigma'(g,x))=\bO_{(g,x)}.
	$$
To prove this, we observe that this amounts to the fact that 
	\begin{quote}
	$\sigma'(g,x)\,(=g.\pi_{\eta}'(x))\in \pi_{\eta}'(X')$ implies $g\in H$.
	\end{quote}
Let us prove this fact.
We write $x=(B_{1},B_{2},i,Z_{1},Z_{2},...,Z_{e})\in \pi_\eta'(X')$.
Recall that the sum of generalized spaces of $B_{1}$ with eigenvalues in the set $Z_{l}$ is $V_{l}^{0}$ for each $l=1,2,...,e$.
The assumption $g.\pi_{\eta}'(x)\in \pi_{\eta}'(X')$ asserts that $g.B_{1}$ also has $V_{l}^{0}$ as the sum of generalized spaces with eigenvalues in $Z_{l}$ for each $l$.
This forces $g\inv(V_{l}^{0})=V_{l}^{0}$ for each $l$, hence $g\in H$.

It remains to prove that $\sigma'$ is submersive onto the image, equivalently the differential $d\sigma'$ at any point has constant rank $\dim G+\dim X'-\dim H$.
This amounts to $\Ker(d\sigma'{}_{(g,x)})=T_{(g,x)}\bO_{(g,x)}$.
By the $G$-equivariancy we may assume that $g$ is the identity element of $G$.
It is direct to observe that both spaces $\Ker(d\sigma'{}_{(g,x)})$ and $T_{(g,x)}\bO_{(g,x)}$ have the following description
    $$
    \Ker (d\sigma'_{(g,x)})
    =\left\{(\xi,\delta x)\in \fg\oplus \bigoplus_{l=1}^e
    \bN_{V^0_l,W} \middle|\xi.x+\delta x=0\right\}=T_{(g,x)}\bO_{(g,x)}.
    $$
This completes the proof.
\qed\vskip.3cm

If $\eta=(k)$, we obtain $X'=Y', X=Y$ and $G=H$. 
In this case, the above lemma gives only the obvious isomorphism.
Otherwise this lemma gives a proper \'etale open subset of $\mu\inv(0)$ associated to $\eta$ arising from the strictly lower instanton numbers.
For the reason, we call this property on $\mu\inv(0)$ \textit{factorization property} and $\pi_{\eta}$ \textit{factorization morphism}.

Since $(S^{\eta}\cc)_{0}$ is an affine variety, both $X,Y$ are also affine schemes.
We denote by $\opi_{\eta}\colon X\git H\to Y\git G$ the induced morphism from $\pi_{\eta}$ between the affine GIT quotients.
Lemma \ref{lem: quotient} (2) says the following:

\begin{cor}\label{cor: factorization isomorphism} 
$\opi_{\eta}$ is an isomorphism.
\end{cor}

We call $\opi_{\eta}$ \textit{factorization isomorphism}.

According to the original argument \cite[Lemma 12.3.2]{MO}, factorization morphism was constructed from the ordinary ADHM data (for $\SU(N)$-instantons).
A similar argument works for $\USp(N/2)$-instantons.
In any cases we denote the factorization morphism and factorization isomorphism by the same symbols $\pi_{\eta},\opi_{\eta}$ respectively.


\subsection{Proof of Theorem \ref{th: main}}
\label{subsec: the proof of main theorem}

We prove the assertions in Theorem \ref{th: main} other than the flatness of $\mu$ by studying scheme structure of the \'etale open subsets $\mu\inv(0)\times_{S^{k}\cc}(S^{\eta}\cc)_{0}$ of $\mu\inv(0)$ for two partitions $\eta=(2^{k/2}),(4,2^{k/2-2})$.
The assertions in the theorem other than normality are deduced from the case $\eta=(2^{k/2})$.
For, the projection image of $\mu\inv(0)\times_{S^{k}\cc}(S^{\eta}\cc)_{0}$ in $\mu\inv(0)$ has the complement $\mu\inv(0)_{\le k/2-1}$.
Hence by Corollary \ref{cor: reduction}, the local geometric properties of $\mu\inv(0)$ outside codimension $1$ locus in the assertions are equivalent to those of $\mu\inv(0)\times_{S^{k}\cc}(S^{\eta}\cc)_{0}$.
Due to Lemma \ref{lem: quotient}, the local geometry of $\mu\inv(0)\times_{S^{k}\cc}(S^{\eta}\cc)_{0}$ now comes from the $k=2$ case via the factorization morphism.
Further in the $k=2$ case (i.e., $\dim V=2$), we have the identification $\mu\inv(0)=\cc^2\times\rho\inv(0)$ where $\rho\colon L(W,V)\to \fg(V),\ i\mapsto ii^*$ since $\fp(V)$ consists of scalars and thus the ADHM equation becomes $[B_1,B_2]=ii^*=0$.

\begin{lem}\label{lem: KP k2N4}
$($\cite[Remark 11.3]{KP}$)$
Suppose $k=2,N\ge4$.
Let $\rho\colon L(W,V)\to\fg(V)$, $i\mapsto ii^{*}$.
Then $\rho\inv(0)$ is a complete intersection and its smooth locus consists of surjective maps.

Moreover the following assertions are true.

$(1)$
If $N=4$, $\rho\inv(0)$ is a reduced scheme with $2$ irreducible components.
These irreducible components are interchanged by the action of any element in $\rO(W)\setminus\SO(W)$.
The smooth locus of $\rho\inv(0)$ is one $\Sp(V)\times\rO(W)$-orbit.

$(2)$
If $N\ge5$, $\rho\inv(0)$ is an irreducible normal variety.
\end{lem}

By this lemma the complement of the codimension $1$ locus $\mu\inv(0)_{\le k/2-1}$ in $\mu\inv(0)$ is a reduced variety.
Therefore the whole scheme $\mu\inv(0)$ itself is also reduced since it is Cohen-Macaulay as a complete intersection (\cite[Prop.\ 5.8.5]{EGA4}).

We prove the irreducibility of $\mu\inv(0)$ for $N\ge5$.
Lemma \ref{lem: KP k2N4} (2) (combined with the factorization property) says $\mu\inv(0)\times_{S^k\cc}(S^\eta\cc)_0$ is irreducible.
Since the projection to $\mu\inv(0)$ is dominant, $\mu\inv(0)$ itself is also irreducible.

We prove the statement of the theorem that the regular locus $\mu\inv(0)^\reg$ is Zariski dense open in $\mu\inv(0)$.
Let $\eta=(2^{k/2})$.
We recall that the factorization morphism $\pi_{\eta}$ composed with the projection to $\mu\inv(0)$ maps the product of the stable-costable quiver representations to stable-costable ones (Lemma \ref{lem: costable}).
Here we used the fact that stability and costability of $\bN$ are equivalent.
Thus we may assume $k=2$.
On the other hand, the surjective maps $i$ in Lemma \ref{lem: KP k2N4} give stable representations in $X$.
Since the costable locus in $\rho\inv(0)$ is Zariski dense open in $\rho\inv(0)$ (Lemma \ref{lem: KP k2N4}), so is $\mu\inv(0)^{\reg}$ in $\mu\inv(0)$.

We describe the irreducible components when $N=4$.
We are still assuming $\eta=(2^{k/2})$.
So $V_{l}^{0}$ are all set to be mutually orthogonal 2-dimensional symplectic vector spaces for $1\le l\le k/2$.
Let $C_{0},C_{1}$ be the two irreducible components of $(\mu^{l})\inv(0)$ coming from those described in Lemma \ref{lem: KP k2N4} (1) for each $1\le l\le k/2$, where $\mu^{l}$ denotes the moment map on $\bN_{V_{l}^{0},W}$.
Let $C_{a_{1}a_{2}\cdots a_{k/2}}$ be the irreducible component in $\mu\inv(0)\times_{S^{k}\cc}(S^{\eta}\cc)_{0}$ corresponding to the product $C_{a_{1}}\times C_{a_{2}}\times\cdots \times C_{a_{k/2}}$ via the factorization morphism where $a_{l}\in\{0,1\}$.
The projection from $\mu\inv(0)\times_{S^{k}\cc}(S^{\eta}\cc)_{0}$ to $\mu\inv(0)$ identifies two components $C_{a_{1}a_{2}\cdots a_{k/2}},C_{b_{1}b_{2}\cdots b_{k/2}}$ if and only if the $k/2$-tuples $(a_{1},a_{2},...,,a_{k/2}),(b_{1},b_{2},...,,b_{k/2})$ coincide after permutation of coordinates.
Therefore the irreducible components of $\mu\inv(0)$ bijectively correspond to the unordered $k/2$-tuples of $0,1$.
This completes the proof of the assertions in the theorem except normality.

To prove normality of $\mu\inv(0)$ for $N\ge5$, we need to look at the case $\eta=(4,2^{k/2-2})$ instead of $\eta=(2^{k/2})$.
The projection image of $\mu\inv(0)\times_{S^k\cc}S^{\eta}\cc$ in $\mu\inv(0)$ has the complement $\mu\inv(0)_{\le k/2-2}$.
Thus by Corollary \ref{cor: reduction}, the local geometric properties of $\mu\inv(0)$ outside codimension $2$ locus in the assertions are equivalent to those of $\mu\inv(0)\times_{S^{k}\cc}(S^{\eta}\cc)_{0}$.
Hence by Serre's criterion, the normality comes from that of $\mu\inv(0)$ in the two cases $k=2,4$.
For $k=2$, this is already done by Lemma \ref{lem: KP k2N4} (2).

\begin{lem}\label{lem: normal k4N5}
If $k=4$ and $N\ge5$, $\mu\inv(0)$ is a normal variety.
\end{lem}
The proof of Lemma \ref{lem: normal k4N5} will appear in \S\ref{sec: normality}.
This finishes the proof of Theorem \ref{th: main} assuming Lemmas \ref{lem: reduction}, \ref{lem: normal k4N5}.


\subsection{Further description on irreducible components for $N=4$ via tensor product}
\label{subsec: further description on irreducible components}

Let $K=\SO(4,\rr)$.
In the previous subsection the irreducible components of $\mu\inv(0)\git G$ are indexed by the $S_{n}$-orbits of $(a_{1},a_{2},...,a_{n})$ where $a_{l}\in\{0,1\}$, where $S_{n}$ denotes the symmetric group of $n$ letters.
Recall from \S\ref{subsec: Geometry of the instanton moduli space} that
the tensor product morphism gives an isomorphism from $\cM^{K}_{(n_{1},n_{2})}$ to each irreducible component of $\mu\inv(0)^{\reg}/\Sp(V)$ for any instanton number $(n_{1},n_{2})$ with $\dim V=2n=2(n_{1}+n_{2})$.
In this subsection we will find the explicit correspondence of the indices $(a_{1},a_{2},...,a_{n})$ and the instanton numbers $(n_{1},n_{2})$ (see Corollary \ref{cor: irreducible components}).
We also express the tensor product morphism in terms of ADHM data (Theorem \ref{th: ADHM tensor product morphism}).

We denote by $T_{(n_{1},n_{2})}\colon \cM^{K}_{(n_{1},n_{2})}\to \mu\inv(0)^{\reg}/\Sp(V)$ the tensor product morphism.
In this section we use the notation $\mu_{k}$ if we need to emphasize $k=\dim V$.
We choose $C_{0}$ among the two irreducible components of $\mu_{2}\inv(0)$ such that $C_{0}^{\reg}/\Sp(\cc^{2})$ coincides with $\Image(T_{(1,0)})$.
Thus $C_{1}$ automatically satisfies $C_{1}^{\reg}/\Sp(\cc^{2})=\Image(T_{(0,1)})$.
Here we need a care: 
$T_{(1,0)}$ maps $\cF\mapsto \cF\otimes\cO^{\oplus2}$ while $T_{(0,1)}$ maps $\cF\mapsto\cO^{\oplus2}\otimes \cF$.
Both $\cF\otimes\cO^{\oplus2}$ and $\cO^{\oplus2}\otimes\cF$ are isomorphic to $\cF^{\oplus2}$ as sheaves, but not as framed sheaves (see Remark \ref{rk: injectivity of tensor product morphism} (2)).

For general $n$, the description of irreducible components follows from the property `tensor product commutes with factorization':
\begin{thm}\label{th: tensor product commutes with factorization}
Let $n=n_{1}+n_{2}$.
There is a commutative diagram
	\begin{equation}\label{eq: factorization diagram}
	\xymatrix{
	\left(\cM^{\USp(1)}_{1}\right)^{n_{1}}\times \left(\cM^{\USp(2)}_{1}\right)^{n_{2}} 
	\ar@{-->}[r]\ar[d]_{(T_{(1,0)})^{n_{1}}\times (T_{(0,1)})^{n_{2}}} & 
	\cM^{\USp(1)}_{n_{1}} \times \cM^{\USp(1)}_{n_{2}}\ar[d]^{T_{(n_{1},n_{2})}}
	\\
	\left(\mu_{2}\inv(0)^{\reg}/\Sp(\cc^{2})\right)^{n_{1}}\times 
	\left(\mu_{2}\inv(0)^{\reg}/\Sp(\cc^{2})\right)^{n_{2}} \ar@{-->}[r] & 
	\mu\inv(0)^{\reg}/\Sp(V)
	}
	\end{equation}
where the horizontal maps are factorization rational maps.
\end{thm}

In fact it is more natural to understand this theorem in terms of the corresponding spaces of maps.
In this formulation our factorization property is reinterpreted into Drinfeld's factorization property.
We will review a relevant theory in \S\ref{app: sec: tensor product commutes with factorization}.
Thus the proof of Theorem \ref{th: tensor product commutes with factorization} is postponed to the last part of \S\ref{app: sec: tensor product commutes with factorization}.

\begin{cor}\label{cor: irreducible components}
The tensor product morphism gives an isomorphim of $\cM^{K}_{(n_1,n_2)}$ with the irreducible component of $\mu\inv(0)^{\reg}/\Sp(V)$ indexed by the unordered $n$-tuple $(a_{1},a_{2},...,a_{n})$ with $n_{1}=\#\left\{l\middle| a_{l}=0\right\}$.
\end{cor}

\proof
By Theorem \ref{th: tensor product commutes with factorization}, the image of the tensor product morphism $T_{(n_{1},n_{2})}$ is the image of $\left(C_{0}^{\reg}/\Sp(\cc^{2})\right)^{n_{1}}\times \left(C_{1}^{\reg}/\Sp(\cc^{2})\right)^{n_{2}}$ via the factorization rational map (the lower horizontal map in the diagram \eqref{eq: factorization diagram}).
According to the index rule in the previous subsection, this image maps via the factorization rational map to the irreducible component indexed by the unordered $n$-tuple in the statement.
\qed\vskip.3cm

We give a quiver description of the tensor product morphism (Theorem \ref{th: ADHM tensor product morphism}).
Let $W_1,W_2$ be $\cc^{2}$ the vector representation of $\USp(1)$.
Let $V_{1},V_{2}$ be orthogonal vector spaces of dimensions $n_{1},n_{2}$ respectively.
Let
    $$
    W:=W_{1}\otimes W_{2},\quad V:=V_{1}\otimes W_{2}\oplus W_{1}\otimes V_{2}.
    $$
Then $V,W$ are symplectic and orthogonal vector spaces with the induced bilinear forms respectively.

Now we have the symplectic vector spaces $\bN_{V_{1},W_{1}}, \bN_{V_{2},W_{2}}, \bN_{V,W}$.
Let us consider the following rational map:
	\begin{equation}\label{eq: ADHM for tensor product}
	\begin{aligned}
	&
	\sT_{(n_{1},n_{2})}\colon \bN_{V_{1},W_{1}}\times \bN_{V_{2},W_{2}}\dashrightarrow \bN_{V,W},
	\quad
	((B_{1,1},B_{1,2},i_{1}),(B_{2,1},B_{2,2},i_{2}))\mapsto
	\\
	&
	\quad
	\left(\begin{pmatrix} B_{1,1}\otimes\Id_{W_{2}} & 0
	\\ 0 & \Id_{W_{1}}\otimes B_{2,1}\end{pmatrix}, 	
	\begin{pmatrix} B_{1,2}\otimes\Id_{W_{2}} & *
	\\ * & \Id_{W_{1}}\otimes B_{2,2}\end{pmatrix},
	\begin{pmatrix}
	i_{1}\otimes\Id_{W_{2}}
	\\
	\Id_{W_{1}}\otimes i_{2}
	\end{pmatrix}
	\right)
	\end{aligned}
	\end{equation}
where $*$ in the above are $n\times n$ matrices are determined as the unique solutions of the equations similar to \eqref{eq: factorization}.
The matrix elements of $*$ are the rational functions in the ones of $B_{l,m}$ and $i_{l}$. 
Hence $\sT_{(n_{1},n_{2})}$ is defined over the locus where the eigenvalue sets of $B_{1,1},B_{1,2}$ are mutually disjoint.
$\sT_{(n_{1},n_{2})}$ commutes with the factorization rational maps as follows:

\begin{lem}\label{lem: sT commutes with factorization}
There is a commutative diagram
	$$ 
	\xymatrix{
	(\bN_{\cc,W_{1}})^{\oplus n_{1}} \times (\bN_{\cc,W_{2}})^{\oplus n_{2}} 
	\ar@{-->}[r] \ar@{-->}[d]_{(\sT_{(1,0)})^{\oplus n_{1}}\times (\sT_{(0,1)})^{\oplus n_{2}}} & 
	\bN_{V_{1},W_{1}}\times \bN_{V_{2},W_{2}}\ar@{-->}[d]^{\sT_{(n_{1},n_{2})}}
	\\
	(\bN_{\cc^{2},W})^{\oplus n_{1}}\times (\bN_{\cc^{2},W})^{\oplus n_{2}}
 	\ar@{-->}[r] & 
	\bN_{V,W}
	}
	$$ 
where the horizontal maps are factorization rational maps.
\end{lem}

\proof
This follows from the definition of the factorization morphisms and $\sT_{(n_{1},n_{2})}$.
We leave the details as an exercise (hint: use the uniqueness of solution of the equation \eqref{eq: factorization}).
\qed\vskip.3cm
Let us regard $\rO(V_{1})\times\rO(V_{2})$ as a subgroup of $\Sp(V)$.
Then $\sT_{(n_{1},n_{2})}$ is $\rO(V_{1})\times\rO(V_{2})$-equivariant.
Thus it induces a rational map from $\mu_{n_{1}}\inv(0)\git\rO(V_{1})\times \mu_{n_{2}}\inv(0)\git\rO(V_{2})$ to $\mu\inv(0)\git\Sp(V)$, where $\mu_{n_{1}},\mu_{n_{2}},\mu$ are the moment maps on $\bN_{V_{1},W_{1}},\bN_{V_{2},W_{2}},\bN_{V,W}$ with respect to $\rO(V_{1}),\rO(V_{2}),\Sp(V)$ respectively.
Lemma \ref{lem: sT commutes with factorization} yields an immediate corollary:

\begin{cor}\label{cor: sT commutes with factorization}
There is a commutative diagram
	\begin{equation}\label{eq: sT factorization diagram}
	\xymatrix{
	(\mu_{1}\inv(0)\git \rO(1))^{n_{1}} \times (\mu_{1}\inv(0)\git\rO(1))^{n_{2}} 
	\ar@{-->}[r] \ar@{-->}[d] 
	& 
	\mu_{n_{1}}\inv(0)\git\rO(V_{1})\times \mu_{n_{2}}\inv(0)\git\rO(V_{2})
	\ar@{-->}[d] 
	\\
	(\mu_{2}\inv(0)\git\Sp(1))^{n_{1}}\times (\mu_{2}\inv(0)\git\Sp(1))^{n_{2}}
 	\ar@{-->}[r] & 
	\mu\inv(0)\git\Sp(V)
	}
	\end{equation}
where the horizontal maps are factorization rational maps and the vertical maps are induced from $\sT_{(1,0)},\sT_{(0,1)},\sT_{(n_{1},n_{2})}$.
\end{cor}

\begin{thm}\label{th: ADHM tensor product morphism}
The tensor product morphism $T_{(n_{1},n_{2})}$ is the induced map from $\sT_{(n_{1},n_{2})}$ on the quotients.
\end{thm}

\proof
The assertion is automatic once we identify the diagrams  \eqref{eq: factorization diagram} and \eqref{eq: sT factorization diagram} over the regular locus.
As both diagrams \eqref{eq: factorization diagram} and \eqref{eq: sT factorization diagram} are commutative and the image of factorization morphism is Zariski open dense, it suffices to check that the induced morphisms from $\sT_{(1,0)}$ and $\sT_{(0,1)}$ are the tensor product morphisms $\cF\mapsto \cF\otimes\cO^{\oplus2}$ and $\cO^{\oplus2}\otimes\cF$ respectively.
Let us check this for $\sT_{(1,0)}$ first.
In this case we have $V_{2}=0$.
Looking at \eqref{eq: ADHM for tensor product}, the image ADHM datum of $\sT_{(1,0)}$ becomes $(B_{1,1}^{\oplus2},B_{1,2}^{\oplus2},i_{1}^{\oplus2})$, which is the ADHM datum of the framed vector bundle $\cF\otimes\cO^{\oplus2}$.
The case for $\sT_{(0,1)}$ is similar.
\qed\vskip.3cm

\begin{rk}\label{rk: injectivity of tensor product morphism}
(1)
In general for any classical groups, the tensor product morphism can be also defined in terms of ADHM data as in \eqref{eq: ADHM for tensor product}.
In our case in this subsection, the tensor product morphism is an isomorphism.
This was shown by the correspondence of instantons in \S\ref{subsec: Geometry of the instanton moduli space}, but \textit{not} by the quiver description.
In the context of spaces of maps, this is also clear. 
See \eqref{eq: obvious identification} (cf.\ the proof of Theorem \ref{th: tensor product commutes with factorization} in \S\ref{app: sec: tensor product commutes with factorization}).

We do not know the purely quiver-theoretic proof of bijectivity of the tensor product morphism.
A difficulty occurs as the quiver-theoretic tensor product morphism is \textit{not} defined over the whole product space.

(2) 
Let us check the dominance of the tensor product morphism purely in terms of quiver representations. 
Due to Theorem \ref{th: tensor product commutes with factorization}, it suffices to consider the case $n=1$.
We only need to prove that the tensor product morphism map from $\cM^{K}_{(0,1)}$ to $C_{1}^{\reg}/\Sp(V)$ because we chose $\cM^{K}_{(1,0)}$ which maps to $C_{0}^{\reg}/\Sp(V)$.
Thus we set $V_{1}=V_{2}=\cc$.
Our claim amounts to the following assertion:
	\begin{quote}
	If $i_{1}\in L(W_{1},V_{1})$ and $i_{2}\in L(W_{2},V_{2})$ are surjective maps, $i_{1}\otimes \Id_{W_{2}}$ and $\Id_{W_{1}}\otimes i_{2}$ are contained in different irreducible components of $\rho\inv(0)$ where $\rho\colon L(W,V)\to \fsp(V),\ i\mapsto ii^{*}$.
	\end{quote}
For, any instanton $\cM^{\USp(1)}_{1}$ has the ADHM datum $(B_{1},B_{2},i)$ such that $B_{1},B_{2}\in\cc$ and $i$ is surjective.
Here we used $ii^{*}=0$ since $\fo(\cc)=0$, so that $i_{1}\otimes \Id_{W_{2}},\Id_{W_{1}}\otimes i_{2}\in\rho\inv(0)$.

Let us prove the assertion.
By Lemma \ref{lem: KP k2N4},  $i_{1}\otimes \Id_{W_{2}},\Id_{W_{1}}\otimes i_{2}$ have the same $\Sp(V)\times\rO(W)$-orbit in $L(W,V)$ since they are surjective maps.
To see that they are contained in different irreducible components, we have to find $\tau\in \rO(W)\setminus\SO(W)$ such that
	$$
	i_{1}\otimes \Id_{W_{2}} = (\Id_{W_{1}}\otimes i_{2}) \tau.
	$$
by Lemma \ref{lem: KP k2N4} (1).

Let $W_{1}=W_{2}=\cc^{2}=\cc\langle e_{1},e_{2}\rangle$ where $e_1,e_2$ are the standard symplectic basis.
We set $i_{1},i_{2}$ to be the same linear map given by $e_{1}\mapsto 1,e_{2}\mapsto 0$.
(Any general surjective map in $L(W_{1},V_{1})$ is contained in the $\Sp(W_{1})$-orbit of $i_{1}$.)
We define $\tau$ as a linear map exchanging the basis elements $e_{1}\otimes e_{1}$ and $e_{2}\otimes e_{2}$.
Then $\tau\in \rO(W)$ but does not preserve the orientation, so that $\tau\in \rO(W)\setminus\SO(W)$.
\end{rk}


\section{Flatness of $\mu$ via modality: proof of Lemma \ref{lem: reduction}}\label{sec: flatness via modality}

The main purpose of the section is to give flatness criteria for moment maps in order to prove Lemma \ref{lem: reduction}.
We need dimension estimate of strata of $\mu\inv(0)$, which will be done in terms of modality in this section.

In \S\ref{subsec: modality}--\S\ref{subsec: flatness}, we study flatness of moment maps under general setting.
We consider two symplectic vector spaces $\cV_{1},\cV_{2}$ with Hamiltonian $G$-actions.
We reduce the flatness problem for $(G,\cV_{1}\oplus\cV_{2})$ to $(G^{x},\cV_{2})$ for $x\in\cV_{1}$ using modality.
In \S\ref{subsec: r2k2}, we turn to our main interest, the $\SO(N)$-data, $N\ge4$.
Recall that the general setting includes $\bN_{V,W}=T^{*}\fp(V)\oplus L(W,V)$, $G=\Sp(V)$.
Thus we will prove the flatness of the moment map for $(\Sp(V)^{x},L(W,V))$.
As we will see in the proof of Lemma \ref{lem: reduction} in \S\ref{subsec: reduction lemma}, the essential part is the flatness for $(\Sp(V)^{x},L(W,V))$ when $x$ is a general nilpotent endomorphism in $\fp(V)$ (see Theorem \ref{th: flat half loop algebra} and Corollary \ref{cor: flatness of mu cVd}).


\subsection{Modality}\label{subsec: modality}

Let $G$ be an algebraic group and $\fg:=\Lie(G)$ the Lie algebra of $G$.
Let $\cV$ be a $G$-scheme.
Let $G.x$ be the $G$-orbit and $G^x:=\{g\in G| g.x=x\}$ the stabilizer subgroup of $x\in \cV$.

We denote by
    $$
    \cV_{(G,s)}:=\{x\in \cV|\, \dim G.x=s\},\quad
    \cV^{(G,s)}:=\{x\in \cV|\, \dim G^x=s\}.
    $$
It is clear that $\cV_{(G,s)}=\cV^{(G,\dim G-s)}$.
By \cite[\S13.1]{EGA4}, the stabilizer dimension is upper-semicontinuous.
Thus $\cV_{(G,s)}$ is a locally closed subvariety of $\cV$ for any $s$.

We define \textit{modality} of $(G,\cV)$ as
    $$
    \modality(G:\cV):=\max_{s\in \zz} (\dim \cV_{(G,s)}-s)
    = \max_{s\in \zz} (\dim \cV^{(G,s)}-\dim G+s)
    $$
(cf.\ \cite[p.198]{PS}).
Here we set dimension of empty set to be $-\infty$.

If $\cV$ is a (linear) representation of $G$, there is also Lie algebra $\fg$-action on $\cV$.
Let $\fg.x$ be the orbit and $\fg^x:=\Lie(G^x)$ the stabilizer Lie algebra where $x\in \cV$.
Then we have $\fg^x=\left\{g\in \fg\middle|\, g.x=0\right\}$.
We denote by $\cV_{(\fg,s)}:=\{x\in \cV|\, \dim \fg.x=s\}$ and $\cV^{(\fg,s)}:=\{x\in \cV|\, \dim \fg^x=s\}$.
Clearly they coincide with $\cV_{(G,s)}$ and $\cV^{(G,s)}$ respectively.


\subsection{Comparison results}

We give basic comparison results of modalities.
The first one is comparison with respect to a subgroup:

\begin{prop}
\label{prop: mod}
Let $H<G$ be an algebraic subgroup. Then $\modality(G:\cV)\le \modality(H:\cV)$.
\end{prop}

\proof
Since $\dim G.x\ge\dim H.x$, we have $\cV_{(G,s)}\subset \bigsqcup_{s'\le s}\cV_{(H,s')}$ for any $s$.
Therefore we obtain
    $$
    \dim \cV_{(G,s)}-s\le \max_{s'\le s} (\cV_{(H,s')}) -s \le \max_{s'\le s}(\dim \cV_{(H,s')} -s')\le \modality(H:\cV) .
    $$
Since $s$ is arbitrary we obtain $\modality(G:\cV)\le\modality(H:\cV)$.  \qed\vskip.3cm

We have another comparison result of modalities of two schemes.
Let $(G_1,\cV_1),(G_2,\cV_2)$ be pairs of an algebraic group and a scheme with some group action.

\begin{prop}\label{prop: sum of modality}
For any $s\ge0$, we have $(\cV_1\times \cV_2)^{(G_1\times G_2,s)}= \bigsqcup_{s=s_1+s_2} \cV_1^{(G_1,s_1)}\times \cV_2^{(G_2,s_2)}$.
Hence,
    $$
    \modality(G_1\times G_2:\cV_1\times \cV_2)=\modality(G_1:\cV_1)+\modality(G_2:\cV_2) .
    $$
\end{prop}

\proof
These are immediate from $(G_1\times G_2).(x_1,x_2)=G_1.x_1\times G_2.x_2$.
\qed\vskip.3cm


\subsection{Flatness of moment maps via modality}
\label{subsec: flatness}

The purpose of the subsection is to give criteria of flatness and irreducibility of moment maps in terms of modality.

Let $\cV$ be a $G$-representation and $\mu_\cV$ be the moment map on $T^*\cV$.
Let $p_1,p_2$ be the first and second projections from $T^*\cV=\cV\oplus\cV\dual$.
We use the abbreviated notations $\cV^{(s)},\cV_{(s)}$ instead of $\cV^{(G,s)},\cV_{(G,s)}$ in this section, as far as comparison with different groups does not appear.

Let us consider the scheme-theoretic inverse $\mu_\cV\inv(\xi)$, $\xi\in \fg\dual$.
It is a scheme cut out by the $\dim G$ equations given by the matrix elements of $\mu_\cV$.

\begin{lem}\label{lem: fibre dim}
For any $x\in \cV$, the restriction of $\mu_\cV$ to $\{x\}\times \cV\dual$ is a linear map dual to the Lie algebra action map $a\colon\fg\to \cV$, $\zeta\mapsto \zeta.x$.

Hence for any given $\xi\in \fg\dual$ and $x\in \cV_{(s)}$, $p_1\inv(x)\cap \mu_\cV\inv(\xi)$ is isomorphic to $\aaa^{\dim\cV-s}$ unless it is empty.
\end{lem}

\proof
The first statement is clear from the definition of moment map.

Therefore $p_1\inv(x)\cap \mu_\cV\inv(0)$ coincides with $\Ker(a\dual)=\Image(a)^\perp$.
So it is isomorphic to $\cc^{\dim\cV-s}$ via $p_2$.
The second statement is proven for $\xi=0$.

If $\xi\neq0$, let us assume $p_1\inv(x)\cap \mu_\cV\inv(\xi)$ is not empty.
We choose any element $(x,y)$ in it.
Then it is isomorphic to an affine subspace $y+\Ker(a\dual)$ in $\cV\dual$ via $p_2$.
\qed\vskip.3cm

The lemma above gives a criterion of flatness of $\mu_\cV$ (\cite[(2.3)]{Pa}):

\begin{prop}
\label{prop: flat criterion}
The following are equivalent.
    \begin{enumerate}
    \item
    $\mu_\cV\inv(0)$ is a complete intersection.
    \item
    $ \dim  \cV_{(s)} -s\le  \dim \cV-\dim G\ \mbox{for all $s$ }. $
    \item[(2')]
    $ \modality(G:\cV)\le  \dim \cV-\dim G.$
    \item
    $ \dim  \cV^{(s)} +s\le  \dim \cV\ \mbox{ for all $s$}. $
    \item
    $\mu_\cV$ is flat.
    \end{enumerate}

Hence if one of the above holds, $\mu_\cV$ is equi-dimensional and Cohen-Macaulay.
\end{prop}

\proof
We prove the equivalence (1)$\Leftrightarrow$(2).
$\mu_\cV\inv(0)$ is a complete intersection if and only if it has dimension $\le 2\dim \cV-\dim G$
if and only if $p_1\inv(\cV_{(s)})\cap \mu_\cV\inv(0)$ has dimension $\le 2\dim \cV-\dim G$ for each $s\ge0$.
By Lemma \ref{lem: fibre dim}, $p_1\inv(\cV_{(s)})\cap \mu_\cV\inv(0)$ has dimension $\dim\cV_{(s)}+\dim\cV-s$.
This proves the equivalence.

The equivalences (2)$\Leftrightarrow$(2') and (2)$\Leftrightarrow$(3) are obvious.

By \cite[Ch.\ III, Exer.\ 10.9]{Hart}, the flatness of $\mu_\cV$ amounts to that every nonempty fibre has dimension $2\dim \cV-\dim G$.
So (4)$\Rightarrow$(1) is immediate.

We prove (2)$\Rightarrow$(4).
We notice that every nonempty fibre $\mu_\cV\inv(\xi)$ has dimension $\ge2\dim\cV-\dim G$ because of the number of defining equations.
So we need to prove its dimension $\le 2\dim\cV-\dim G$.
By Lemma \ref{lem: fibre dim}, every nonempty fibre of the restriction of $p_1$ to $p_1\inv(\cV_{(s)})\cap \mu_\cV\inv(\xi)$ has dimension $\dim\cV-s$ for each $s\ge0$.
Thus $p_1\inv(\cV_{(s)})\cap \mu_\cV\inv(\xi)$ has dimension $\le \dim\cV_{(s)}+\dim\cV-s$.
Now the assumption (2) asserts its dimension $\le 2\dim\cV-\dim G$.
We are done.
\qed\vskip.3cm


\subsection{Modality of $L(\cc^{r},\cc^{k})$ with respect to the current Lie algebra of $\fsl_{2}$}
\label{subsec: r2k2}

Let $T:=\cc^{2}$ with the standard symplectic form.
Let $V_{n}=T^{n}=\cc^{k}$.
Let $r\ge2$ and $\cV_{n}:=L(\cc^{r},V_{n})$.

First we study the case when $n=1$.
This case will be the initial inductive step for $n\ge1$.
Let $G:=\Sp(V_{1})=\SL_{2}$.

Let us consider a stratification of $\cV_{1}$ by rank:
    $$
    \cV_{1}^l:=\left\{f\in \cV_{1}\middle| \rank f=l\right\}
    $$
where $l=0,1,2$.
This stratification coincides with the stratification $\cV_{0}^{(s)}$ as follows:
    $$
    \cV_{1}^{(3)}=\cV_{1}^0=0,\quad \cV_{1}^{(1)}=\cV_{1}^1,\quad
    \cV_{1}^{(0)}=\cV_{1}^2.
    $$

\begin{prop}\label{prop: k2r2}
$(1)$
All the strata $\cV_{1}^{(3)},\cV_{1}^{(1)},\cV_{1}^{(0)}$ are irreducible varieties of dimension $0,r+1,2r$ respectively.

Hence $\modality(G:\cV_{1})=2r-3$ and $\mu_\cV$ is flat.

$(2)$
If $r=2$, the modality is attained by $\cV_{1}^{(1)}$ and $\cV_{1}^{(0)}$

$(3)$
If $r\ge3$, the modality is attained by $\cV_{1}^{(0)}$.
\end{prop}

\proof
For each $0\le l\le 2$, we have a surjective morphism $\cV_{1}^l\to \Gr(l,2)$, $f\mapsto \Image(f)$.
It is a Zariski locally trivial fibration over the irreducible variety $\Gr(l,2)$, whose fibre at $\Image(f)$ is the locus of surjective maps in $L(\cc^r,\Image(f))$.
Therefore $\cV_1^l$ is an irreducible variety of dimension $rl+2(l-2)$.
Now all the statements are clear from the the flatness criterion (Proposition \ref{prop: flat criterion}).
\qed\vskip.3cm

Now we study the case when $n\ge2$.
We define natural current Lie algebra $\fsl_{2}[z]$-actions on $V_n$ and $\cV_n$ as follows.
The $\gl_{2}$-action on the vector representation $T$ induces the $\gl_{2}[z]$-action on $T\otimes\cc[z]$.
We call $T\otimes\cc[z]$ the \textit{vector representation} of $\gl_{2}[z]$.
Let $\cc[z]/(z^{n})$ be the truncated polynomial algebra.
We identify $V_{n}=T^{\oplus n}=T\otimes \cc[z]/(z^{n})$ by regarding the $m^{\mathrm{th}}$ direct summand $T$ as $T\otimes z^{m-1}$.
Note that $\cc[z]$ is the center of $\gl_{2}[z]$.
Thus $z^{n}\gl_{2}[z]$ is the two-sided ideal in $\gl_{2}[z]$ generated by $z^{n}$.
The quotient algebra $\gl_{2}[z]/z^{n}\gl_{2}[z]$ naturally acts on $V_{n}$.
Note that the center $\cc[z]$ acts also on the Lie subalgebra $\fsl_{2}[z]$ and that $z^{n}\fsl_{2}[z]$ is a Lie algebra two-sided ideal.
Let $\fg_{n}:=\fsl_{2}[z]/z^{n}\fsl_{2}[z]$.
Then we have the induced Lie algebra $\fg_{n}$-action on $V_{n}$.
There are the induced $\gl_{2}[z]$- and $\fsl_{2}[z]$-actions on $L(\cc^{r},T\otimes\cc[z])$ and thus the $\fg_{n}$-action on $\cV_{n}=L(\cc^{r},T\otimes\cc[z]/(z^{n}))$.  

The main theorem of this subsection is as follows:
\begin{thm}\label{th: flat half loop algebra}
$\modality(\fg_{n}:\cV_{n})=(2r-3)n$.
\end{thm}

The proof of the theorem will appear at the end of the subsection after preliminary steps.
This theorem combined with Proposition \ref{prop: flat criterion} yields an immediate corollary:

\begin{cor}\label{cor: flatness of mu cVd}
The moment map $\mu_{\cV_{d}}\colon T^{*}\cV_{n}\to (\fg_{n})\dual$ is flat for any $d\ge0$.
\end{cor}

In fact, $\fg_{n}$ will be realized as the Lie algebra $\Lie(\Sp(V_{n})^{x})$ for a generic nilpotent endomorphism $x\in\fp(V_{n})$ as we will see in Lemma \ref{lem: fg}.

\begin{rk}\label{rk: current algebra orbit}
Note that any $\fsl_2[z]$-orbits in $V_{n}$ and $\cV_{n}$ coincide with the $\fg_{n}$-orbits respectively since $z^{n}$ annihilates any element.
Thus Theorem \ref{th: flat half loop algebra} assures $\modality(\fsl_{2}[z]:\cV_{n})=\modality(\fg_{n}:\cV_{n})=(2r-3)n$.
\end{rk}

Since there are vector space identifications  
	$$
	V_{n}=T\otimes\cc[z]/(z^{n}),\quad 
	\cV_{n}=L(\cc^{r},T)\otimes\cc[z]/(z^{n}),\quad \fg_{n}=\fsl_{2}\otimes\cc[z]/(z^{n}),
	$$
the elements of $V_{n},\cV_{n},\fg_{n}$ are written as polynomials with coefficients in $T,L(\cc^{r},T),\fsl_{2}$ respectively.
So we use the notation $x_{m}=\Coeff_{m}(x)$ for the coefficient of $z^m$ in a given polynomial $x$.

We stratify $\cV_{n}$ as follows:
	$$
	\cV_{n}^{l}:=\left\{x\in \cV_{n}\middle|\rank (x_{0})=l\right\}
	$$
where $l=0,1,2$.
This is an immediate generalization of the stratification $\cV_1^l$ of $\cV_{1}$ in the $n=1$ case.

We consider the evaluation map   
$\cV_n\otimes \cc^{r}\to V_n$.
This map intertwines the $\gl_{2}[z]$-action on $\cV_{n}$, so there is a natural map $\gl_2[z]\otimes\cV_n\otimes \cc^r \to V_n$.
The image of $\xi\otimes x\otimes w\in V_{n}$ via this map will be denoted by $\xi x w$, where $\xi\in \gl_{2}[z]$, $x\in \cV_{n}$ and $w\in\cc^{r}$.
We also use the notation $\xi x\in\cV_{n}$ and $x w\in V_{n}$ for the natural maps $\gl_{2}[z]\otimes \cV_{n}\to\cV_{n}$ and $\cV_{n}\otimes \cc^{r}\to V_{n}$ respectively.

For $n\ge m$, $\cc[z]/(z^m)$ is naturally a $\cc[z]/(z^n)$-module.
Both the multiplication map
	$$
	z^{n-m}\colon \cc[z]/(z^m)\to \cc[z]/(z^n),\quad f\mod(z^m)\mapsto z^{n-m}f\mod(z^n)
	$$
and the truncation map
	$$
	\cc[z]/(z^n)\to \cc[z]/(z^m),\quad f\mod(z^n)\mapsto f\mod(z^m)
	$$
are $\cc[z]$-module homomorphisms.
These induce the multiplication maps by $z$ and the truncation maps:
	\begin{equation*}
	\begin{aligned}
	&
	z\colon \fg_{n-1}\to \fg_n,\quad z\colon \cV_{n-1}\to \cV_{n}
	\\
	&
	\tau_{\fg_n}\colon \fg_n\to \fg_{n-1},\quad \tau_{\cV_n}\colon \cV_n\to \cV_{n-1}
	\end{aligned}
	\end{equation*}
which are $\fsl_2[z]$- or $\gl_2[z]$-module homomorphisms appropriately.
From polynomial expression, the following identifications as vector spaces are clear:
    $$
    \Ker (\tau_{\fg_n})=\fsl_2\otimes z^{n-1},
    \quad
    \Ker (\tau_{\cV_n})=L(\cc^r,T)\otimes z^{n-1},
    \quad
    \tau_{\cV_n}(\cV_n^l)=\cV_{n-1}^l
    $$
for any $l\ge0$.

We define the \textit{minimal degree} of $x\in\fg_n$ and $\cV_n$ as follows:
    $$
    \mindeg x:=\max\left\{0\le m\le n\middle| z^{n-m}x=0\right\}.
    $$
It is nothing but the smallest exponent $m$ of $z^m$ with nonzero coefficient in the polynomial expression of $x$, if $x\neq0$.
Note that $\mindeg x=n$ if and only if $x=0$.

We use several vector space identifications of $\Ker(z)$ in $\cV_n$:
    $$
    \Ker(z)=\Ker(\tau_{\cV_{n}})=\{x\in \cV_{n}| \mindeg x\ge n-1\} =z^{n-1}\cV_n=z^{n-1}\cV_1.
    $$
This is immediate from the polynomial expressions.
We have also a similar identification on $\fg_n$.

We now describe the stabilizer $\fg_n^x$ or the orbit $\fg_n x$ for $x\in\cV_n^l$ where $l=0,1,2$.
We start from $l=2$ and then $l=0,1$.

\begin{prop}\label{prop: rk2}
$\fg_n^x=0$ for any $x\in \cV_n^2$.
\end{prop}

\proof
Let $\xi\in \fg_n\setminus0$.
Let $m:=\mindeg \xi$.
Since $x_0\colon \cc^r\to T$ is surjective, the composite $\xi_mx_0$ and thus $\xi x$ are both nonzero.
This means that no nonzero $\xi$ annihilates $x$.
\qed\vskip.3cm

\begin{prop}\label{prop: rk0}
$(\cV_{n}^0)_{(\fg_n,s)}=z (\cV_{n-1})_{(\fg_{n-1},s)}$ for any $s\ge0$.
\end{prop}

\proof
The multiplication map $z\colon \cV_{n-1}\to \cV_{n}^{0}$ is an injective $\fsl_2[z]$-module homomorphism.
Hence it gives an isomorphism between the $\fsl_2[z]$-orbits (Remark \ref{rk: current algebra orbit}), which proves the proposition.
\qed\vskip.3cm

The last case to study is $l=1$.
It is more complicated than $l=0,2$.

Let $x\in \cV_n^1$.
We describe $\fg_n^x$ first.
Since $x_0\colon \cc^r\to T$ is of rank 1, there is a nonzero element $t\in \fsl_2$ such that $\fsl_2^{x_0}=\cc t$.
Thus $\fg_n^x$ is not trivial since $t\otimes z^{n-1}$ annihilates $x$.
In fact $\fg_d^x$ is generated by one element as a $\cc[z]$-module due to the lemma below.
We denote by
	$$
	m_x:=\min_{\xi\in \fg_n^x} (\mindeg \xi).
	$$

\begin{lem}
\label{lem: principal gen}
Let $x\in \cV_n^1$ and $\xi\in \fg_n^x$ with $\mindeg \xi=m_x$.
Then we have $\fg_n^x=\cc[z]\xi$.
In particular $\dim \fg_n^x=n-m_x$.
\end{lem}

\proof
Let $\xi'\in \fg_n^x$.
It suffices to check $\xi'\in \cc[z]\xi$, because the opposite inclusion $\cc[z]\xi\subset \fg_{n}^{x}$ is obvious.
Let $m':=\mindeg \xi'$.
If $m'=n+1$ then $\xi'=0$ and thus $\xi'\in \cc[z]\xi$.
We use the induction on $m'$; the case $m'=n+1$ is the initial induction hypothesis and we assume the lemma is true for $m'+1,m'+2,...,n+1$. 

Since both coefficients $\xi_{m_x},\xi'_{m'}$ annihilate $x_0$, there exists $c\in \cc^*$ such that $\xi'_{m'}=c\xi_{m_x}$.
Here we used the above observation that $\dim \fsl_2^{x_0}=1$.
Now $\xi'-cz^{m'-m_x}\xi$ has minimal degree strictly larger than $m'$ and annihilates $x$.
By the induction hypothesis, it is contained in $\cc[z]\xi$, which proves $\fg_{n}^{x}\subset \cc[z]\xi$.
This completes the proof.
\qed\vskip.3cm

Let $x\in \cV_n^1$ as before.
By $\xi_x$, we denote an element $\xi$ in the above lemma.
The choice of $\xi_x$ is unique up to multiplication by units in $\cc[z]/(z^{n})$ according to the lemma.

If $m_x\ge1$ then $\xi_x$ is divisible by $z$ in $\fg_n$.
We choose any $\xi_x'\in \fg_n$ with $z\xi_x'=\xi_x$.
The choice of $\xi_x'$ is unique up to addition by elements of $z^{n-1}\fg_{n-1}$.

\begin{prop}\label{prop: rk1}
Let $n\ge2$ and $x\in \cV_{n}^1$.
Then the stabilizer Lie algebra of $\tau_{\cV_{n}}(x)$ in $\fg_{n-1}$ is identified as
    $$
    \fg_{n-1}^{\tau_{\cV_n}(x)}=\left\{
        \begin{array}{lll}
        \cc[z]\tau_{\fg_n}(\xi_x) && \mbox{if $m_x=0$}
        \\
        \cc[z]\tau_{\fg_n}(\xi_x') && \mbox{if $m_x\ge1$.}
        \end{array}
    \right.
    $$

Hence if $s=n$, then $\tau_{\cV_n}((\cV_n^1)^{(\fg_n,s)})\subset (\cV_{n-1}^1)^{(\fg_n,s-1)}$.
If $s<n$, then $\tau_{\cV_n}((\cV_n^1)^{(\fg_n,s)})\subset (\cV_{n-1}^1)^{(\fg_n,s)}$.
\end{prop}

\proof
Let $\fg':=\fg_{n-1}^{\tau_{\cV_n}(x)}$ for short during the proof.

Since $\tau_{\cV_n}(x)\in \cV_{n-1}^1$, there is $\xi\in \fg_{n-1}$ such that
$\fg'=\cc[z]\xi$ by Lemma \ref{lem: principal gen}.
We need to show that such $\xi$ can be chosen as in the statement.
By Lemma \ref{lem: principal gen}, this amounts to showing that $\tau_{\fg_n}(\xi_x)$ and $\tau_{\fg_n}(\xi_x')$ give the smallest minimal degree of elements of $\fg'$ in the respective cases.

Truncation of both sides of the identity $\xi_xx=0$ gives $\tau_{\fg_n}(\xi_x)\tau_{\cV_n}(x)=\tau_{\cV_{n}}(\xi_{x}x)=0$, so $\tau_{\fg_n}(\xi_x)\in \fg'$.
If $m_x=0$ then $\tau_{\fg_n}(\xi_x)$ has minimal degree 0.
So it generates $\fg'$ as a $\cc[z]$-module by Lemma \ref{lem: principal gen}.
If $m_x\ge1$, it is direct to check that $\tau_{\fg_n}(\xi_x')$ also annihilates $\tau_{\cV_n}(x)$ and has minimal degree $m_x-1$.
To conclude that it generates $\fg'$ as a $\cc[z]$-module, it suffices to check that $\mindeg \xi'\ge m_x-1$ for all nonzero $\xi'\in \fg'$ due to Lemma \ref{lem: principal gen}.
Since $z\xi'$ annihilates $x$ in $\cV_d$, we have $\mindeg z\xi' \ge m_x$.
Therefore, $\mindeg z\xi'=1+\mindeg \xi'\ge m_{x}$.

The assertions of the inclusions in the statement follow from the above description of $\fg'$ and the dimension formula $\dim \cc[z]\xi_x=n-m_x$ in  Lemma \ref{lem: principal gen}.
\qed\vskip.3cm

\begin{lem}\label{lem: fibre dim of trunc}
If $n\ge2$, all the fibres of the restriction of $\tau_{\cV_n}$ to
$(\cV_n^1)^{(\fg_n,n)}$ have dimension $r+1$.
\end{lem}

\proof
We fix any $x\in (\cV_n^1)^{(\fg_n,n)}$ and a $\cc[z]$-module generator $\xi_x$ of $\fg_d^x$.
Then for any $y\in \tau_{\cV_n}\inv\tau_{\cV_n}(x)$, the following are equivalent:
	\begin{enumerate}
	\item
	$y\in (\cV_n^1)^{(\fg_n,n)}$.
	\item
	There exists $\xi_y\in \fg_n^y$ of minimal degree 0.
	\item
	There exists  $\xi_y\in \fg_n^y$ of minimal degree 0 such that $\xi_x\equiv \xi_y\mod (z^{n-1})$.
	\end{enumerate}
The implication $(3)\Rightarrow(2)$ is obvious.
The equivalence $(1)\Leftrightarrow (2)$ is clear from Lemma \ref{lem: principal gen}.
We prove the implication $(1)\Rightarrow (3)$.
By the observation after Lemma \ref{lem: principal gen}, there is $f(z)\in \cc[z]$ with $f(0)\neq0$ such that $f(z)\xi_y\equiv\xi_x\mod (z^{n-1})$, because $\xi_xy\equiv0\mod (z^{n-1})$ as $x\equiv y\mod (z^{n-1})$.
Thus by replacing $\xi_y$ into $f(z)\xi_y$, we are done.

Now we can describe the fibre $\tau_{\cV_n}\inv\tau_{\cV_n}(x)\cap (\cV_n^1)^{(\fg_n,n)}$ using the above equivalences.
We use the following polynomial expressions:
 	\begin{equation*}
	\begin{aligned}
	&
	x=\sum_{m=0}^{n-1} x_m\otimes z^m,\quad
	\xi_x=\sum_{m=0}^{n-1}\xi_{x,m}\otimes z^m,
	\\
	&
	y=y_{n-1}\otimes z^{n-1}+\sum_{m=0}^{n-2}x_m\otimes z^m,
	\quad
	\xi=\xi_{n-1}\otimes z^{n-1}+ \sum_{m=0}^{n-2}\xi_{x,m}\otimes z^m.
	\end{aligned}
	\end{equation*}
By the equivalence $(1)\Leftrightarrow (3)$, the above fibre $\tau_{\cV_n}\inv\tau_{\cV_n}(x)\cap (\cV_n^1)^{(\fg_n,n)}$ is isomorphic to the locus of $y_{n-1}\in L(\cc^r,T)$ such that there exists $\xi_{n-1}\in \fsl_2$ satisfying $\xi y=0$.
Note that
	$$
	\xi y=\xi y-\xi_xx=\left(\xi_{x,0}(y_{n-1}-x_{n-1})+(\xi_{n-1}-\xi_{x,n-1})x_0\right)\otimes z^{n-1}.
	$$
Thus our problem is reduced to the general question on the first order deformation of linear maps $A\colon \bV_2\to \bV_3$ and $B\colon \bV_1\to \bV_2$ with constraint $AB=0$, where $\bV_1,\bV_2,\bV_3$ are finite dimensional vector spaces.
We are looking for the first order deformations $\delta A,\delta B$ of $A,B$ respectively such that
	\begin{equation}\label{eq: AB}
	A(\delta B)+(\delta A) B=0.
	\end{equation}
We observe that a pair $(\delta A,\delta B)$ is a solution of \eqref{eq: AB} if and only if $\delta B(\Ker B)\subset \Ker A$.
So the second factor $\delta B$ of the solutions $(\delta A,\delta B)$ forms a vector space of dimension $(\dim \Ker A)(\dim \Ker B)+(\dim \bV_2)(\dim \bV_1-\dim \Ker B)$. Here the second summand comes from the choice of linear maps from a complementary subspace of $\Ker B$ to $\bV_2$.
In our case, $A=\xi_{x,0},B=x_0$, $\bV_{1}=\cc^{r}$, $\bV_{2}=\bV_{3}=T$.
Thus the dimension formula of the vector space of $\delta B$ is computed as $1\cdot(r-1)+2\cdot1=r+1$.
This completes the proof.
\qed\vskip.3cm

Now we are ready to prove Theorem \ref{th: flat half loop algebra}.
\vskip.3cm

\textit{Proof of Theorem \ref{th: flat half loop algebra}.}
The equality in the theorem amounts to the following inequality
	\begin{equation}\label{eq: dim cVdl}
	\dim(\cV_n^{l})^{(\fg_n,s)}+s\le 2rn
	\end{equation}
for all $s\ge0$ and $l=0,1,2$.
We notice that Proposition \ref{prop: flat criterion} assures that the equality holds for at least one pair $(s,l)$.

Let $l=2$ first.
By Proposition \ref{prop: rk2}, the equality holds for $s=0$ and $(\cV_n^{l})^{(\fg_n,s)}=\emptyset$ for the other values $s$.

For $l=0,1$, we use the induction on $n$, so we assume that the inequality holds up to $n-1$.
The initial induction step $n=1$ was proven in Proposition \ref{prop: k2r2}.

Suppose $l=1$.
There are the two subcases: $s<n$ and $s=n$.
We assume $s<n$ first.
Since the truncation map $\tau_{\cV_{n}}\colon\cV_{n}\to\cV_{n-1}$ is a linear map with kernel dimension $2r$ and $\tau_{\cV_{n}}\left((\cV_n^l)^{(\fg_n,s)}\right)\subset (\cV_{n-1}^l)^{(\fg_{n-1},s)}$ (Proposition \ref{prop: rk1}), we obtain the dimension estimate
$\dim (\cV_n^l)^{(\fg_n,s)}\le \dim (\cV_{n-1}^l)^{(\fg_{n-1},s)}+2r$.
We add $s$ to both sides of this inequality.
As a result we get 
 	$$
	\dim (\cV_n^l)^{(\fg_n,s)}+s\le \dim (\cV_{n-1}^l)^{(\fg_{n-1},s)}+s+2r\le 2r(n-1)+2r=2rn,
	$$ 
where we used the induction hypothesis to the second inequality.
We get the inequality \eqref{eq: dim cVdl} for $l=1$ and $s<n$.

We assume $s=n$ secondly.
By Proposition \ref{prop: rk1}, $\tau_{\cV_{n}}\left((\cV_n^l)^{(\fg_n,s)}\right)\subset (\cV_{n-1}^l)^{(\fg_{n-1},s-1)}$.
By Lemma \ref{lem: fibre dim of trunc}, any fibre of the restriction of $\tau_{\cV_{n}}$ to $(\cV_n^l)^{(\fg_n,s)}$ has dimension $r+1\,(\le 2r-1)$.
Thus a similar argument gives $\dim (\cV_n^l)^{(\fg_n,s)}+s\le \dim (\cV_{n-1}^l)^{(\fg_{n-1},s-1)}+(s-1)+2r\le 2r(n-1)+2r$.
This completes the proof for $l=1$.

If $l=0$, using Proposition \ref{prop: rk0}, we rewrite \eqref{eq: dim cVdl} into the inequality 
$\dim \cV_{n-1}^{(\fg_{n-1},s-3)}+s\le 2rn$.
The left hand side is less than or equal to $2r(n-1)+3$ by the induction hypothesis.
Since $2r(n-1)+3< 2rn$, we obtain \eqref{eq: dim cVdl}.
\qed\vskip.3cm


\subsection{Proof of Lemma \ref{lem: reduction}}\label{subsec: reduction lemma}

We need to show that every nonempty fibre of $\tmu$ has dimension $\le \dim \bN-\fg(V)-n$ since the target space $\fg(V)\times S^{n}\cc$ is smooth (\cite[Ch.\ III, Exer.\ 10.9]{Hart}).
Since there is the contraction $\cc^{*}$-action on $\bN$, we consider only the zero-fibre by the method of associated cones (\cite[II.4.2]{Kr}).
So we need to check $\mu\inv(0)\cap p\inv(\fp(V)^{\nilp})$ has the above dimension $\dim \bN-\fg(V)-n$, where $p\colon \bN\to\fp(V)$ is the first projection and $\fp(V)^{\nilp}$ denotes the locus of nilpotent endomorphisms in $\fp(V)$.

Note that $\fp(V)^{\nilp}$ is a finite union of $\Sp(V)$-orbits (Lemma \ref{lem: G}).
By the dimension estimate of the fibres of the projection $p_1\colon \mu\inv(0)\cap p\inv(\Sp(V).B)\to \Sp(V).B\times L(\cc^r,V)$ (Lemma \ref{lem: fibre dim}), we need to check
	\begin{equation}\label{eq: modality for nilp}
	\modality(\Sp(V):\Sp(V).B\times L(\cc^{r},V))\le\dim L(\cc^{r},V)-3n
	\end{equation}
for any nilpotent endomorphism $B\in \fp(V)$ where $r\ge2$.

On the other hand there is a reciprocal equality:
	\begin{equation}\label{eq: mod mod}
	\modality(\Sp(V):\Sp(V).B\times L(\cc^{r},V))=\modality(\Sp(V)^{B}:L(\cc^{r},V)).
	\end{equation}
This is deduced as follows:
Let $(\Sp(V).B\times L(\cc^{r},V))^{(\Sp(V),s)}\to \Sp(V).B$ be the projection where $s\ge0$.
The fibre of this map at $B$ is $L(\cc^{r},V)^{(\Sp(V)^{B},s)}$ and the other fibres are all isomorphic to this.
Thus $\dim (\Sp(V).B\times L(\cc^{r},V))^{(\Sp(V),s)}=\dim L(\cc^{r},V)^{(\Sp(V)^{B},s)}+\dim \Sp(V).B$.
To both sides of this equality we add $s-\dim \Sp(V)=s-\dim \Sp(V)^{B}-\dim\Sp(V).B$.
Then  we obtain \eqref{eq: mod mod}.

By Lemma \ref{lem: G}, for any $B\in\fp(V)^{\nilp}$, $V$ is decomposed into symplectic subspaces $V_{l}$ such that $B|_{V_{l}}$ is a generic nilpotent endomorphism in $\fp(V_{l})$.
Then we have
    \begin{equation}\label{eq: two comparison}
    \begin{aligned}
	\modality(\Sp(V)^B:L(\cc^r,V))
	&
	\le
    \modality\left(\prod_{l=1}^n\Sp(V_l)^{B|_{V_{l}}}:\bigoplus_{l=1}^nL(\cc^r,V_l)
    \right)
    \\
    &
    =\sum_{l=1}^n \modality\left(\Sp(V_l)^{B|_{V_{l}}}:L(\cc^r,V_l)\right)
    \end{aligned}
    \end{equation}
where we used the two comparison results (Propositions \ref{prop: mod}, \ref{prop: sum of modality}).
Therefore by \eqref{eq: modality for nilp} and \eqref{eq: mod mod}, it suffices to show $\modality(\Sp(V)^{B}:L(\cc^{r},V))=(2r-3)n$ for a generic nilpotent endomorphism $B\in\fp(V)$.

The above equality is immediate from Theorem \ref{th: flat half loop algebra} under suitable identification of $(V,B)$ as follows:
 Let $G:=\Sp(V)^{B}$ for short.
By Lemma \ref{lem: G}, we can identify $(V,B)$ and $(T\otimes \cc[z]/(z^{n}),z)$, where
the symplectic form on the right hand side is given as
	$$
	(\of_1,\of_2)=\Res\left(\frac{(f_1,f_2)_T}{z^{n}}\right).
	$$
Here $f_1,f_2\in T\otimes \cc[z]$, $\of_1,\of_2$ are their truncations and the residue is taken at $z=0$.
The pairing $(\,,\,)_T$ is the $\cc[z]$-bilinear extension of the symplectic form of $T$ over $T\otimes \cc[z]$.
It is clear that the above residue form depends only on the classes $\of_1,\of_2$ but not on $f_1,f_2$ themselves.

\begin{lem}\label{lem: fg}
Under the above identification of $(V,B)$, the Lie algebra $\fg=\Lie(\Sp(V)^{B})$ coincides with the truncated current algebra $\fsl_2[z]/z^{n}\fsl_{2}[z]$.
\end{lem}

\proof
For any $\xi\in \fsl_2[z]/(z^{n})$ and $\of_1,\of_2\in V$, we have
	$$
	(\xi \of_1,\of_2)=\Res\left(\frac{(\xi \of_1,\of_2)_T}{z^{n}}\right)
	=-\Res\left(\frac{(\of_1,\xi\of_2)_T}{z^{n}}\right)=-(\of_1,\xi \of_2).
	$$
Thus $\fsl_2[z]/z^{n}\fsl_{2}[z]$ is a subset of $\fg$.

We show the opposite inclusion.
Since any endomorphism $\xi$ of $V$ in $\fg$ commutes with $z$, it is an element of $\gl_2[z]/z^{n}\gl_{2}[z]$.
Since $\xi$ satisfies the above equalities for any pairs $\of_1,\of_2\in V$, it always comes from $\fsl_2[z]/z^{n}\fsl_{2}[z]$. 
\qed \vskip.3cm

This finishes the proof of Lemma \ref{lem: reduction}.

\begin{rk}\label{rk: SO3}
Lemma \ref{lem: reduction} does \textit{not} hold for $N=3$ if $k\ge4$.
For, by the factorization property it is enough to consider the case $k=4$.
By \cite[Theorem 4.3]{Choy}, $\mu\inv(0)$ is a complete intersection with two irreducible components.
By \cite[Corollary 8.9]{Choy}, one irreducible component is contained in $\mu\inv(0)\cap (E_kp)\inv(\Delta S^{k/2}\cc)$ where $p\colon \bN_{V,W}\to \fp(V)$ is the first projection.
Since $\Delta S^{k/2}\cc$ has codimension $1$ in $S^{k/2}\cc$, the restriction of $E_kp$ to $\mu\inv(0)$ is not equidimensional, which shows $\tmu$ is not flat.
\end{rk}


\section{Normality of $\mu\inv(0)$: proof of Lemma \ref{lem: normal k4N5}}
\label{sec: normality}

Let $V,W$ be a symplectic and orthogonal vector spaces of dimension $4,5$ respectively as before.
In this case we prove normality of $\mu\inv(0)$ (Lemma \ref{lem: normal k4N5}).
For $\dim W>5$, the normality of $\mu\inv(0)$ follows from the $\dim W=5$ case by the base change argument as was noticed in \S\ref{subsec: moment maps}.
See the proof of this argument in \S\ref{App: normality}.


\subsection{Proof of normality in Lemma \ref{lem: normal k4N5}}\label{subsec: proof of normality}

In order to prove the normality of $\mu\inv(0)$, we will check that the regular locus $\mu\inv(0)^{\reg}$ has the complement of codimension $\ge2$ (use Serre's criterion as $\mu\inv(0)$ is Cohen-Macaulay).
First we consider the locus of $(B_{1},B_{2},i)$ such that $ii^{*}$ is not nilpotent endomorphism.
This locus will turn out to be a Zariski open subset of $\mu\inv(0)^{\reg}$.
We will see also that the regular element outside this locus is Zariski dense open in the locus. 

We will use the Kraft-Procesi theory on nilpotent pairs (\cite{KP1}\cite{KP}).
A pair $(i,i^{*})$ is called \textit{nilpotent pair} if $ii^{*}$ is a nilpotent endomorphism.

We start the details from some basic properties on $\fp(V), L(W,V)$.
Let $\fp:=\fp(V),\fg:=\fg(V)$ for short.
For any $(B_{1},B_{2})\in\fp^{\oplus2}$, the square of commutator $([B_{1},B_{2}])^{2}$ is always a scalar endomorphism by a direct calculation. 
So for $(B_{1},B_{2},i)\in\mu\inv(0)$, we have $(ii^{*})^{2}$ is also a scalar.
Let $V_{a},W_{a}$ be the generalized $a$-eigenspaces of $ii^{*},i^{*}i$ respectively, where $a\in\cc$.
Then we have $i(W_{a})\subset V_{a},\ i^{*}(V_{a})\subset W_{a}$ (\cite[Lemma 5.5]{Choy}).
So by the fact that $(ii^{*})^{2}$ is scalar, the set of eigenvalues of $ii^{*}$ is $\{\pm a\}$ for some $a\in \cc$.

Suppose first that $a\neq0$.
Then both $i,i^{*}$ are isomorphisms between $V_{a}$ and $W_{a}$ (\cite[Lemma 5.5]{Choy}).
Note that $V_{a},V_{-a}$ (resp.\ $W_{a},W_{-a}$) are isotropic subspaces dual to each other via $(\,,\,)_{V}$ (resp.\ $(\,,\,)_{W}$) (\cite[Lemma 5.2]{Choy}).
So $V_{a},V_{-a},W_{a},W_{-a}$ are all 2-dimensional.
Thus the orthogonal complement $(W_{a}\oplus W_{-a})^\perp$ of $W_{a}\oplus W_{-a}$ in $W$ is automatically $W_{0}$.
For, $(W_{a}\oplus W_{-a})^\perp$ is an eigenspace of $i^{*}i$ as it is 1-dimensional and $i^{*}i$-invariant.
Let $b\in\cc\setminus\{\pm a\}$ be the eigenvalue of the restriction of $i^{*}i$ to $(W_{a}\oplus W_{-a})^\perp$.
Then $b$ should be $0$ because otherwise $i$ is isomorphism between $W_{b},V_{b}$, which cannot happen as there is no eigenvalue $b$ of $ii^{*}$.

We describe first 
	$$
	Z:=\{i\in L(W,V)|\mbox{$(ii^{*})^{2}$ is nonzero scalar}\}.
	$$

\begin{lem}\label{lem: description of Z}
$\dim Z=17$.
\end{lem}

\proof
We consider a morphism $Z\to\cc^{*},\ i\mapsto (ii^{*})^{2}$.
We denote by $Z_{a}$ the fibre at $a\in \cc^{*}$.
We fix any $a\in\cc^{*}$ and $i\in Z_{a^{2}}$.
We know that $i$ gives rise to the $2$-dimensional isotropic generalized eigenspaces $V_{a},V_{-a},W_{a},W_{-a}$ and the isomorphism $i\colon L(W_{a},V_{a})$.
We observe further that
	\begin{enumerate}
	\item
	the isomorphism $i|_{W_{a}}\colon W_{a}\to V_{a}$ is dual to $i^{*}|_{V_{-a}}\colon V_{-a}\to W_{-a}$  	 under the identifications $V_{a}\dual\cong V_{-a},W_{a}\dual\cong W_{-a}$ via $(\,,\,)_{V},(\,,\,)_{W}$ 	 respectively,
	\item
    $(i|_{W_{a}})(i^{*}|_{V_{a}})= a\Id_{V_{a}}$,
	\item
	$(i|_{W_{-a}})(i^{*}|_{V_{-a}})= -a\Id_{V_{a}}$.
	\end{enumerate}
Let $(\Gr^{\iso}(2,V)\times\Gr^{\iso}(2,V))_{0}$ be the Zariski open subset of the product $\Gr^{\iso}(2,V)\times\Gr^{\iso}(2,V)$ consisting of the pairs $(V_{1},V_{2})$ of mutually complementary $2$-dimensional isotropic subspaces in $V$ dual to each other via $(\,,\,)_{V}$.
We define $(\Gr^{\iso}(2,W)\times\Gr^{\iso}(2,W))_{0}$ similarly.
We define a morphism
 	$$
	\begin{aligned}	
	&
	Z_{a^{2}}\to (\Gr^{\iso}(2,V)\times \Gr^{\iso}(2,V))_{0}\times (\Gr^{\iso}(2,W)\times\Gr^{\iso}(2,W))_{0},
	\\
	&
	i\mapsto (V_{a},V_{-a},W_{a},W_{-a}).
	\end{aligned}
	$$
Conversely given $(V_a,V_{-a},W_a,W_{-a})$, one can recover $i$ once we specify an arbitrary isomorphism $i|_{W_a}\colon W_a\to V_a$.
Thus we have
    $$
    \dim Z_{a^{2}}=\dim L(V_{a},W_{a})+2(\dim \Gr^{\iso}(2,V)+\dim \Gr^{\iso}(2,W))=16,
    $$
where we used $\dim \Gr^{\iso}(2,V)=4,\ \dim\Gr^{\iso}(2,W)=2$.

Since the morphism $Z\to\cc^{*}$ has the fibre dimension $16$, $Z$ has dimension $17$.
\qed\vskip.3cm

Secondly we describe
	$$
	Z':=\{i\in L(W,V)|(ii^{*})^{2}=0\},
	$$
which is the complement of $Z$.
We apply Kraft-Procesi's theory on nilpotent pairs via $ab$-diagrams \cite{KP}.
We use the notations $a,b$ as basis elements of $W,V$ respectively.
According to the theory (\cite[Theorem 6.5]{KP}), there is a bijection between the set of $\Sp(V)\times\rO(V)$-orbits of $i$ in $L(W,V)$ with nilpotent $ii^{*}$ and the set of $ab$-diagrams satisfying some rules.
We list all the possible $ab$-diagrams of $i$ with $(ii^{*})^{2}=0$ in Table \ref{table: possible ab}.
	\begin{table}[h]%
	\footnotesize
	\caption{$ab$-diagrams, $\Delta_{ab}$, dimensions of $\Sp(V).ii^{*},\rO(W).i^{*}i,(\Sp(V)\times\rO(W)).i$}
	\label{table: possible ab}\centering %
	\begin{tabular}{llllll}
	$ab$-diagram & $\Delta_{ab}$ & $\dim\Sp(V).ii^{*}$ & $\dim \rO(W).i^{*}i$ & $\dim(\Sp(V)\times\rO(W)).i$
 	\vspace{0.1cm} \\ \hline \vspace{0.1cm}
	${\xy <1cm,0cm>:
   	(1,0.4)*{ababa},
    	(1,0)*{bab},
    	(1,-0.4)*{a}
    	\endxy}$
	& 0 & 6 & 6 & 16
	\vspace{0.1cm} \\ \hline \vspace{0.1cm}
	${\xy <1cm,0cm>:
   	(1,0.4)*{abab},
    	(1,0)*{baba},
    	(1,-0.4)*{a}
    	\endxy}$
	& 0 & 6 & 4 & 15
	\vspace{0.1cm} \\ \hline \vspace{0.1cm}
	${\xy <1cm,0cm>:
   	(1,0.4)*{bab},
    	(1,0)*{bab},
    	(1,-0.4)*{a^{3}}
    	\endxy}$
	& 0 & 6 & 0 & 13
	\vspace{0.1cm} \\ \hline \vspace{0.1cm}
	${\xy <1cm,0cm>:
   	(1,0.4)*{ababa},
    	(1,0)*{ab},
    	(1,-0.4)*{ba}
    	\endxy}$
	& 0 & 4 & 6 & 15
	\vspace{0.1cm} \\ \hline \vspace{0.1cm}
 	${\xy <1cm,0cm>:
   	(1,0.4)*{ababa},
    	(1,0)*{b^{2}},
    	(1,-0.4)*{a^{2}}
    	\endxy}$
	& 4 & 4 & 6 & 13
	\vspace{0.1cm} \\ \hline \vspace{0.1cm}
	${\xy <1cm,0cm>:
   	(1,0.4)*{bab},
    	(1,0)*{aba},
    	(1,-0.4)*{aba}
    	\endxy}$
	& 2 & 4 & 4 & 13
	\vspace{0.1cm} \\ \hline \vspace{0.1cm}
	${\xy <1cm,0cm>:
   	(1,0.4)*{bab},
    	(1,0)*{ab},
    	(1,-0.4)*{ba},
	(1,-0.8)*{a^{2}}
    	\endxy}$
	& 0 & 4 & 0 & 12
	\vspace{0.1cm} \\ \hline \vspace{0.1cm}
	${\xy <1cm,0cm>:
   	(1,0.4)*{bab},
    	(1,0)*{b^{2}},
    	(1,-0.4)*{a^{4}}
    	\endxy}$
	& 6 & 4 & 0 & 9
	\vspace{0.1cm} \\ \hline \vspace{0.1cm}
	${\xy <1cm,0cm>:
   	(1,0.4)*{aba},
    	(1,0)*{aba},
    	(1,-0.4)*{b^{2}},
	(1,-0.8)*{a}
    	\endxy}$
	& 2 & 0 & 4 & 11
	\vspace{0.1cm} \\ \hline \vspace{0.1cm}
	${\xy <1cm,0cm>:
   	(1,0.4)*{ab},
    	(1,0)*{ba},
    	(1,-0.4)*{b^{2}},
	(1,-0.8)*{a^{3}}
    	\endxy}$
	& 6 & 0 & 0 & 7
	\vspace{0.1cm} \\ \hline \vspace{0.1cm}
	${\xy <1cm,0cm>:
   	(1,0.4)*{b^{4}},
    	(1,0)*{a^{5}}
    	\endxy}$
	& 20 & 0 & 0 & 0
	\vspace{0.1cm} \\ \hline \vspace{0.1cm}
	\end{tabular}
	\end{table}

In the table we used the abbreviated notations in $ab$-diagrams, for instance 
	$$
	a^{3}={\xy <1cm,0cm>:
   	(1,0.4)*{a},
    	(1,0)*{a},
    	(1,-0.4)*{a}
    	\endxy}$$
According to the $ab$-diagrams in the table, $Z'$ is the disjoint union of eleven $\Sp(V)\times\rO(W)$-orbits.
We also computed dimensions of these orbits in the table by using the formula \cite[Proposition 7.1]{KP}:
	$$
	\dim(\Sp(V)\times \rO(W)).i= \frac12(\dim \Sp(V).ii^{*}+\dim \rO(W).i^{*}i+ \dim V.\dim W-\Delta_{ab}).
	$$

Let $m\colon\fp^{\oplus2}\to\fg,\ (B_{1},B_{2})\mapsto [B_{1},B_{2}]$ (the commutator map).
We claim that dimension of the nonempty $m$-fibre of any nonzero element (resp.\ $m\inv(0)$) is $5$ (resp.\ $8$).
This claim is proven in the explanation in \cite[(8.6)]{Choy}, but we sketch the proof as we will use a similar argument in the proof of Proposition \ref{prop: k2}.
Let $\fp'$ be the space of trace-free endomorphisms in $\fp$.
We have decomposition $\fp=\cc\oplus\fp'$ using splitting by the linear embedding of the scalar endomorphisms $\cc$ in $\fp$ and the trace map on $\fp$.
We denote by $m'$ the restriction of $m$ to $\fp'{}^{\oplus2}$.
By the universality of the exterior product, $m'$ factors through the wedge product map $\fp'{}^{\oplus2}\to\Lambda^{2}\fp'$, $(B_{1},B_{2})\mapsto B_{1}\wedge B_{2}$.
In fact the factored map $\Lambda^{2}\fp'\to \fg$ is an isomorphism.
See \cite[(8.6)]{Choy} for verification of this fact. 
If $B_1\wedge B_2\neq0$, it defines a $2$-dimensional subspace of $\fp'$ equipped with a volume form.
Therefore $m'{}\inv(B_1\wedge B_2)$ is $3$-dimensional.
On the other hand $B_1\wedge B_2=0$ means that $B_1$ and $B_2$ are proportional. 
Therefore $m'{}\inv(0)$ is $6$-dimensional.
Since $m$ is the pull-back of $m'$ via the projection $\fp^{\oplus2}=\cc^{2}\oplus\fp'{}^{\oplus2}\to\fp'{}^{\oplus2}$, the claim is proven.

We decompose $\mu\inv(0)$ into two parts $\mu\inv(0)\cap p\inv(Z')$ and $\mu\inv(0)\cap p\inv(Z)$ where $p$ is the projection from $\bN$ to $L(W,V)$.
We further decompose the former locus using the decomposition of $Z'$ into the eleven $\Sp(V)\times\rO(W)$-orbits as above.
We denote these orbits by $O_{1},O_{2},...,O_{11}$ from the top in order in Table \ref{table: possible ab}.
Note that for any $i\in L(W,V)$, $\mu\inv(0)\cap p\inv(i)\cong m\inv(ii^{*})$.
Note also that the Young diagram corresponding to $ii^{*}$ is obtained by deleting $a$ in the $ab$-diagram of $i$. 
Thus $\mu\inv(0)\cap p\inv(O_{l})$ has dimension $\dim O_{l}$ plus 8 or 5 depending on $ii^{*}=0$ or not for each $1\le l\le 11$.
So $\dim\mu\inv(0)\cap p\inv(O_{l})=21,20,18,20,18,18,17,14,19,15,8$ for $l=1,2,...,11$ respectively.
Similarly the locus $\mu\inv(0)\cap p\inv(Z)$ has dimension $\dim Z+5=22$.

Now we show every element of $\mu\inv(0)\cap p\inv(Z)$ and $\mu\inv(0)\cap p\inv(O_{1})$ is regular.
First any element $(B_{1},B_{2},i)\in \mu\inv(0)\cap p\inv(Z)$ is costable because $i^{*}$ is injective.
Since stability and costability are equivalent for SO-data, it is regular.
Secondly we also prove $(B_{1},B_{2},i)\in\mu\inv(0)\cap p\inv(O_{1})$ is costable.
We take $B_{1},B_{2},i$ explicitly as follows:
Let $\{e_{1},e_{2},e_{3},e_{4}\},\{f_{1},f_{2},...,f_{5}\}$ be bases of $V,W$ whose pairing matrices $(e_{l},e_{m})_{V}$ and $(f_{l},f_{m})_{W}$ are
	$$
	\begin{pmatrix}
		0 & 1 & 0 & 0 \\ -1 & 0 & 0 & 0 \\ 0 & 0 & 0 & 1 \\ 0 & 0 & -1 & 0
	\end{pmatrix},
	\quad
	\begin{pmatrix}
		0 & 1 & 0 & 0 & 0 \\ 1 & 0 & 0 & 0 & 0 \\ 0 & 0 & 1 & 0 & 0 \\ 0 & 0 & 0 & 0 & 1 \\ 0 & 0 & 0 & 1 & 0
	\end{pmatrix}
	$$
respectively.
We define
	$$
	B_1:=\begin{pmatrix}
		\frac12 & 0 & 0 & 0 \\ 0 & \frac12 & 0 & -\frac12 \\ -\frac12 & 0 & -\frac12 & 0 \\ 0 & 0 & 0 & -\frac12 			 \end{pmatrix},
	\quad
	B_2:= \begin{pmatrix}
		0 & 0 & 0 & 0 \\ 0 & 0 & 0 & 1 \\ 1 & 0 & 0 & 0 \\ 0 & 0 & 0 & 0
	\end{pmatrix}
	$$
in $\fp$ with respect to $e_{1},e_{2},e_{3},e_{4}$.
We define
	$$
	i:=\begin{pmatrix} 0 & 0 & 0 & 1 & 0 \\ 0 & -1 & 0 & 0 & 0 \\ 1 & 0 & 0 & 0 & 0 \\ 0 & 0 & 0 & 1 & 0
	\end{pmatrix}
	$$
with respect to bases $e_{1},e_{2},e_{3},e_{4},f_{1},f_{2},...,f_{5}$.
Then with respect to the same bases, we have
	$$
	i^{*}=\begin{pmatrix}
		1 & 0 & 0 & 0 \\ 0 & 0 & 0 & 1 \\ 0 & 0 & 0 & 0 \\ 0 & 0 & 0 & 0 \\ 0 & 1 & -1 & 0
	\end{pmatrix}.
	$$
It is direct to check that $[B_{1},B_{2}]+ii^{*}=0$.
Also by ranks of the alternating composites $ii^{*}ii^{*}...$ and $i^{*}ii^{*}i...$, the $ab$-diagram of $i$ is the first one of Table \ref{table: possible ab} and thus $i\in O_{1}$.

We check the above triple $(B_{1},B_{2},i)$ is costable, hence regular.
From the above explicit form of $i^{*}$, we have $\Ker(i^{*})=\cc\langle e_{2}+e_{3}\rangle$.
However from $B_{1}(e_{2}+e_{3})=\frac12(e_{2}-e_{3})$, there is no nonzero $B_{1},B_{2}$-invariant subspace in $\Ker(i^{*})$.
Therefore $(B_{1},B_{2},i)$ is costable.

Let $\SL_2'$ be the subgroup of affine linear transforms on $\cc^2$ generated by $\SL_2$ and the translations.
There is a natural $\SL_2'$-action on $\fp^{\oplus2}$ by
    $$
    \begin{aligned}
    &
    \begin{pmatrix} a & b \\ c& d \end{pmatrix}.(B_1,B_2):=(aB_1+bB_2,cB_1+dB_2),
    \\
    &
    (a,b).(B_1,B_2):=(B_1+a\Id_V,B_2+b\Id_V).
    \end{aligned}
    $$
We already proved in the above that any fibre of $p\colon \mu\inv(0)\cap p\inv(O_{11})\to O_{11}$
is an $\SL_2'$-orbit, because it is isomorphic to $m\inv(ii^{*})$ and $ii^{*}\neq0$.
Since $O_{11}$ is an $\Sp(V)\times\rO(W)$-orbit, 
$\mu\inv(0)\cap p\inv(O_{11})$ is an $\SL_2'\times\Sp(V)\times\rO(W)$-orbit
(see \cite[proof of Corollary 8.9]{Choy}).
Therefore every element of it is also regular.

This completes the proof of normality in Lemma \ref{lem: normal k4N5}.
\qed\vskip.3cm


\section{Variants of Lemma \ref{lem: reduction} for various ADHM data}
\label{sec: reduction lemma for ADHM}

In this section we prove a similar statement with Lemma \ref{lem: reduction} for ordinary ADHM data and Sp-data.
As an application we study geometry of $\mu\inv(0)$ in \S\ref{subsec: remark}.



\subsection{Lemma \ref{lem: reduction} for the ordinary ADHM data and Sp-data}
\label{subsec: reduction lemma for ADHM data}

For the ordinary ADHM data (resp.\ Sp-data), we fix vector spaces (resp.\ an orthogonal and symplectic vector spaces) $V,W$ of dimensions $k\ge0,N\ge1$ (resp.\ $k\ge0, N\in 2\zz_{\ge0}$).
Let $E_{k}\colon\gl(V)\to S^{k}\cc$, $B\mapsto$ the set of eigenvalues of $B$ counted with multiplicities.

\begin{lem}\label{lem: reduction for ADHM}
The morphism
	$$
	\bM_{V,W}\to \gl(V)\times S^{k}\cc, \quad x=(B_{1},B_{2},i,j)\mapsto (\mu(x),E_k(B_{1}))
	$$
is flat.

If $V,W$ are orthogonal and symplectic vector spaces, the morphism $\bN_{V,W}\to \fg(V)\times S^{k}\cc$ obtained by the restriction is also flat.
\end{lem}

The first assertion for $N=1$ is proven by Gan-Ginzburg \cite[Proposition 2.3.2]{GG} with a different method.

\proof
The idea of proof goes as in \S\ref{subsec: reduction lemma}.
By smoothness of the target spaces and the method of associated cones, we need to show the dimensions of the zero fibres of the above morphisms are $\dim \bM_{V,W}-k^{2}-k$ and $\dim \bN_{V,W}-k(k-1)/2-k$ respectively.
Due to the reciprocal equality as in \eqref{eq: mod mod}, this amounts to the following inequalities on the modalities
	\begin{equation*}
    \begin{aligned}
    &
	\modality(\GL(V)^{B}:L(W,V))\le\dim L(W,V)-k,
    \\
    &
    \modality(\rO(V)^{B}:L(W_{L},V))\le\dim L(W_{L},V)
    \end{aligned}
    \end{equation*}
for any nilpotent endomorphism $B$ in $\gl(V),\fp(V)$ respectively, where $W_{L}$ is a maximal isotropic subspace of $W$.

We prove the above inequality for usual ADHM data first.
We set first $B$ to be regular nilpotent.
Then $(V,B)$ is identified with $(\cc[z]/(z^{k}),z)$.
Thus $\gl(V)^{B}=\cc[z]/(z^{k})$.
The stratum of elements with $\gl(V)^{B}$-orbit dimension $\le s$ is given as follows:
   $$
    \left(W\dual\otimes\cc[z]/(z^k)\right)_{(\le s)}=\left\{
        \begin{array}{lll}W\dual\otimes z^{k-s} \cc[z]/(z^k) && \mbox{if $s\le k$}
        \\
        0 && \mbox{if $s>k$.}
        \end{array}
    \right.
    $$
Since the above stratum has dimension $sN$ for each $s\le k$, we obtain
    $$
    \modality\left(\cc[z]/(z^k):W\dual\otimes\cc[z]/(z^k)\right) = kN-k.
    $$
If $B$ is not regular, we decompose $V$ into subspaces on which $B$ gives a regular endomorphism by Jordan normal form.
By the two comparison results as in \eqref{eq: two comparison}, we also have $\modality(\GL(V)^{B}:L(W,V))\le kN-k$.
This finishes the proof for usual ADHM data.

Next we prove the case of Sp-data.
We set first $B$ to be generic nilpotent.
Then $(V,B)$ and the orthogonal form on $V$ are identified with $(\cc[z]/(z^{k}),z)$ and $(\overline{f},\overline{g})_{V}=\Res(fg/z^{k})$ respectively (Lemma \ref{lem: G}).
See \S\ref{subsec: reduction lemma} for the notation here. 
Since $\gl(V)^{B}=\cc[z]/(z^{k})$ consists of only symmetric endomorphisms, we have $\fg(V)^{B}=0$ and thus the modality inequality is automatic.
If $B$ is not generic, we decompose $V$ into mutually orthogonal subspaces on which $B$ gives a generic endomorphism (Lemma \ref{lem: G}).
This finishes the proof for the Sp-data.
\qed\vskip.3cm


\subsection{Geometry of $\mu\inv(0)$}
\label{subsec: remark}

Using Lemma \ref{lem: reduction for ADHM} and the factorization property, we can describe local geometric structures of $\mu\inv(0)$ as was done in \S\ref{subsec: factorization}.
We denote by $\mu_{k}$ the moment map for the instanton number $k$ as in \eqref{eq: stratification of Uhlenbeck space}.
For the usual ADHM data with rank $1$ and instanton number $1$, $\mu_1\inv(0)$ is $\cc^2$ times the union of two distinct lines in $\cc^2$.
From this, $\mu_k\inv(0)$ for rank $1$ and arbitrary instanton number $k$, is a reduced scheme with $k+1$ irreducible components.
This fact was proven in \cite[Theorem 1.1.2]{GG} with a different approach.
From the stratification of the Uhlenbeck space \eqref{eq: stratification of Uhlenbeck space}, only one irreducible component of $\mu_k\inv(0)$ contains the regular locus.

For the case rank $N\ge2$, $\mu_1\inv(0)$ is $\cc^2$ times the minimal nilpotent $\GL(N)$-orbit closure $\ocO$ in $\gl(N)$.
Since $\ocO$ is an irreducible reduced variety of dimension $2N-1$ with the cone singularity, $\mu_k\inv(0)$ in the case $N\ge2$ and any $k$, is also irreducible and reduced.
For the case rank $N\ge2$ and instanton number $2$, the factors $(B_1,B_2)$ of $\mu_2\inv(0)$ are no more commuting pairs, so we need the geometric description of $\mu_2\inv(0)$ as in the proof of Lemma \ref{lem: normal k4N5}.
We can apply the same argument for Sp-data.

\begin{prop}\label{prop: k2}
$(1)$ For the usual ADHM data with $N\ge2$ and any $k$, $\mu_k\inv(0)$ is an irreducible normal variety.

$(2)$ For the Sp-data with $N\ge0$ and $k=2$, $\mu_2\inv(0)$ is the cone over a smooth irreducible variety.
Hence $\mu_k\inv(0)$ is also an irreducible normal variety for any $k$.
\end{prop}

\proof
It suffices to consider the case $k=2$ due to the base change argument.
Let $V:=\cc^2$ the standard $2$-dimensional orthogonal vector space.
Let $\SL_2'$ be the group generated by $\SL_2$ and the translations as in \S\ref{subsec: proof of normality}.
We have commuting diagrams:
    $$
    \xymatrix{
        & \Lambda^2\fsl_2 \ar[d]^\cong
        \\
        (\fsl_2)^{\oplus2} \ar[r]_{m_{\fsl_2}}\ar[ur]^{m'_{\fsl_2}}  & \fsl_2
        }\quad\quad\quad
    \xymatrix{
        & \Lambda^2\fp' \ar[d]^\cong
        \\
        (\fp')^{\oplus2} \ar[r]_{m_{\fp'}}\ar[ur]^{m'_{\fp'}} & \fg
        }
    $$
where $m_\bullet$ is the commutator map, $m'_\bullet$ is the wedge product map, $\fp'$ is the trace-free part of $\fp(V)$, $\fg:=\fg(V)$.
Hence we can use a similar argument in \S\ref{subsec: proof of normality}.

For the usual ADHM data, we suppose $N=2$.
The case $N\ge3$ comes from the base change argument (\S\ref{app: sec: tensor product commutes with factorization}).
There is the (Zariski dense) open stratum in $\mu_2\inv(0)$ from $(\mu_1\inv(0)\times\mu_1\inv(0))\times_{S^2\cc}(S^\eta\cc)_0$ via the factorization morphism where $\eta=(1,1)$.
We denote this open stratum by $S$.
We already know that $S$ and thus $g.S$ are irreducible and normal for any $g\in \SL_{2}'$ due to the explanation above the proposition.
Hence $\SL_{2}'.S=\bigcup_{g\in \SL(2)'}g.S$ is also irreducible and normal.

Now the first item for usual ADHM data comes the claim: $\mu_2\inv(0)\setminus\SL_{2}'.S$ has codimension $2$.
We prove the claim.
Let $(B_1,B_2,i,j)\in\mu_2\inv(0)\setminus\SL_{2}'.S$.
First we check that further if $B_1,B_2\in\fsl_{2}$, they are linearly dependent nilpotent endomorphisms.
Since the linear dependency of the pairs in $\fsl_{2}^{\oplus2}$ is stable under the simultaneous conjugation, we may assume that $B_{1}=X$ where we denote the Chevalley basis by
	$$
	H:=\begin{pmatrix} 1 & 0 	\\ 0 & -1 \end{pmatrix},\quad X:=\begin{pmatrix} 0 & 1 	\\ 0 & 0 \end{pmatrix},\quad Y:=\begin{pmatrix} 0 & 0 	\\ 1 & 0 \end{pmatrix}.
	$$
Here we used the assumption that $B_{1}$ has only one eigenvalue.
We need to show $B_{2}$ is a scalar multiple of $X$.
Using the explicit form $B_{2}=aH+bX+cY$, unless $a=c=0$, it is direct to check that there is $g\in \SL(2)$ such that the first factor of $g.(B_{1},B_{2})$ has two distinct eigenvalues.
This proves the linear dependency and nilpotency.
By the linear dependency, for $(B_1,B_2,i,j)\in\mu_2\inv(0)\setminus\SL_{2}'.S$, we obtain $ij=0$.
Due to Kraft-Procesi's theory of nilpotent pairs (\cite[Proposition 5.3]{KP1}), the pairs $(i,j)$ satisfying $ij=0$ form finitely many $\GL(V)\times\GL(W)$-orbits and among those orbits the following $ab$-diagram
    $$
    {\xy <1cm,0cm>:
   	(1,0.4)*{aba},
    	(1,0)*{b}
    	\endxy}
    $$
gives the unique orbit with the largest dimension $5$.
Hence the locus $\mu_2\inv(0)\setminus\SL_{2}'.S$ has dimension $10$ because the pairs $(B_1,B_2)$ of linearly dependent nilpotent endomorphisms form a $3$-dimensional variety.
Since $\dim\mu_2\inv(0)=12$, the proof of claim is done.
This finishes the case of usual ADHM data.

The second item for the Sp-data follows by looking at the zero-fibre of $m'_{\fp'}$.
It is the cone over $U^{\oplus2}$ where $U$ denotes the universal bundle over the Grassmannian $\Gr(1,2)$.
We omit the detail.
\qed\vskip.3cm

\begin{rk}\label{rk: two proofs}
We gave two ways of proof for normality of $\mu\inv(0)$ for various ADHM data in the instanton number $2$ case so far: (1) detailed analysis of the complement of the open stratum $S$ using Kraft-Procesi's theory as in \S\ref{subsec: proof of normality}, (2) complement codimension $\ge2$ of the locus $\SL_{2}'.S$ as in the above.
\end{rk}



\appendix


\section{Normality of $\mu\inv(0)$ via base change argument}\label{App: normality}

Let $\cV_{1},\cV_{2}$ are  $G$-representations.
Let $\mu_{1},\mu_{2}$ be the moment maps on $\cV_{1},\cV_{2}$ respectively.
Recall that $\mu=\mu_{1}+\mu_{2}$.
We consider the Hamiltonian triples $(\cV_{1},G,\mu_{1}), (\cV_{1}\oplus \cV_{2},G,\mu)$ as in \S\ref{subsec: moment maps}.

Suppose that $\mu_{1}\inv(0)$ is a normal scheme with the expected dimension $\dim \cV_{1}-\dim G$.
In this section, we will deduce that $\mu\inv(0)$ is an irreducible normal variety.
By the assumption, $\mu\inv(0)$ has also the expected dimension $\dim \cV_{1}+\dim\cV_{2}-\dim G$.
Thus it is a Cohen-Macaulay scheme as a complete intersection.
Next we notice that both $\mu\inv(0),\mu_{1}\inv(0)$ are cones, thus connected spaces.
Thus once $\mu\inv(0)$ is a normal scheme, it is automatically an irreducible normal variety.

We start to check that $\mu\inv(0)$ is a normal scheme.
By Serre's criterion, we need to check that the smooth locus $\mu\inv(0)^{\sm}$ has the complement of codimension $\ge2$.
Note that $\mu\inv(0)^{\sm}$ is the locus of $x$ such that the differential $d\mu_{x}$ is surjective because $T_x(\mu\inv(0))=\Ker(d\mu_{x})$.

On the other hand, $\mu_{1}\inv(\xi)$ is an irreducible normal variety for any $\xi\in\fg\dual$, since so is the central fibre $\mu_{1}\inv(0)$.
This fact follows from the method of associated cones \cite[II.4.2]{Kr}.
(The setting in \cite[II.4.2]{Kr} is different from ours, so a modification is needed. 
See \cite[Theorem D.1]{Choy_Thesis}.) 
Therefore $\mu_{1}\inv(\xi)^{\sm}=\{y\in \cV_{1}|\Coker((d\mu_{1})_{y})=0\}$ has the complement of codimension $\ge2$ in $\mu_{1}\inv(\xi)$.

Let $p\colon \mu\inv(0)\to \cV_{2}$ be the projection to the second factor.
Then $p\inv(z)=\mu_{1}\inv(\mu_{2}(z))$.
Note that $\Image(d\mu)_{(y,z)}$ contains $\Image(d\mu_{1})_{y}$ for any $z\in \cV_{2}$.
This implies that the locus $\left\{(y,z)\in\cV_{1}\oplus\cV_{2}\middle| y\in p\inv(z)^{\sm}\right\}$ is a sublocus of $\mu\inv(0)^{\sm}$.
Since for each $z\in \cV_{2}$, $p\inv(z)^{\sm}=\mu_{1}\inv(\mu_{2}(z))^{\sm}$ has the complement of codimension $\ge2$ in $p\inv(z)$, the union $\bigcup_{z\in\cV_{2}}p\inv(z)^{\sm}$ also has the complement of codimension $\ge2$ in $\mu\inv(0)$.
Therefore $\mu\inv(0)\setminus\mu\inv(0)^{\sm}$ has dimension $\ge2$.
As a result, $\mu\inv(0)$ is an irreducible normal variety.


\section{Tensor product commutes with factorization: proof of Theorem \ref{th: tensor product commutes with factorization}}
\label{app: sec: tensor product commutes with factorization}

We prove Theorem \ref{th: tensor product commutes with factorization} by identifying the instanton spaces with spaces of maps.
The proof will appear at the end of this section.
We use the foundation of the theory of spaces of maps in Braverman-Finkelberg-Gaitsgory's work \cite{BFG}.

Let us consider a stack $\cT$  over $\cc$ with an open substack $\cT^{0}$ such that the complement $\cT\setminus \cT^{0}$ is a Cartier divisor of $\cT$.
We suppose that $\cT^{0}$ is isomorphic to the one-point scheme $\Spec\cc$.
Let $\Map(\pp^{1},\cT)$ be the functor of based maps from $\pp^{1}$ to $\cT$.
This parametrizes the maps $f\colon\pp^{1}\to \cT$ sending $\infty\mapsto \cT^{0}$, where $\infty$ denotes the point at infinity $\pp^{1}\setminus\cc$.
Unless $f$ is a constant map, the Cartier divisor $f\inv(\cT\setminus\cT^{0})$ in $\cc=\pp^{1}\setminus\infty$ defines a point of the symmetric product of $\cc$.
Thus we obtain a morphism from the substack of non-constant maps to the symmetric product $S^{\ge1}\cc=\coprod_{\ge1}S^{n}\cc$.
Let $\Map^{n}(\pp^{1},\cT)$ be the inverse image via this morphism of $S^{n}\cc$ where $n\ge1$.
We denote by $\Map^{0}(\pp^{1},\cT):=\Map(\pp^{1},\cT)\setminus\Map^{\ge1}(\pp^{1},\cT)=\Map(\pp^{1},\cT^{0})$, which is a one-point scheme. 
We set $S^{0}\cc=\Spec\cc$. 
Thus we obtain a morphism
	$$
	\sE_{\cT}^{n}\colon \Map^{n}(\pp^{1},\cT)\to S^{n}\cc.
	$$
Here when $n=0$, $\sE^{0}_{\cT}$ is set to be the constant map.
Taking the union over all $n\ge0$, we obtain $\sE_{\cT}\colon \Map(\pp^{1},\cT)\to S^{\bullet}\cc$, where $S^{\bullet}\cc=\coprod_{n\ge0}S^{n}\cc$.

The following is called Drinfeld's factorization property \cite[Proposition 2.17]{BFG}:
\begin{prop}\label{prop: Drinfeld's factorization property}
For each $n_{1},n_{2}\ge0$, there is a natural isomorphism 
	$$
	\begin{aligned}	
	\varpi_{\cT}^{n_{1},n_{2}}\colon
	&
	(\Map^{n_{1}}(\pp^{1},\cT)\times \Map^{n_{2}}(\pp^{1},\cT))_{0}
	\cong \Map^{n_{1}+n_{2}}(\pp^{1},\cT)\times_{S^{n_{1}+n_{2}}\cc}
	(S^{n_{1}}\cc\times S^{n_{2}}\cc)_{0},
	\\
	& (f,g)\mapsto (f\cup g,f\inv(\cT\setminus\cT^{0}),g\inv(\cT\setminus\cT^{0})).
	\end{aligned}
	$$
\end{prop}

Here the notation are given as follows: 
If $n_{1},n_{2}\ge1$, $(S^{n_{1}}\cc\times S^{n_{2}}\cc)_{0}$ denotes the locus of divisors with disjoint supports (\S\ref{subsec: factorization}) and 
	$$
	\begin{aligned}
	&
	(\Map^{n_{1}}(\pp^{1},\cT)\times \Map^{n_{2}}(\pp^{1},\cT))_{0}
	\\
	&
	:=(\Map^{n_{1}}(\pp^{1},\cT)\times \Map^{n_{2}}(\pp^{1},\cT))
	\times_{S^{n_{1}}\cc\times S^{n_{2}}\cc}(S^{n_{1}}\cc\times S^{n_{2}}\cc)_{0}.
	\end{aligned}
	$$
The map $f\cup g$ denotes the pasting of $f,g$ defined as follows: 
\begin{equation*}
 f\cup g= \left\{
 	\begin{array}{lll} f && \mbox{outside $g\inv(\cT\setminus\cT^{0})$}
	\\
	g && \mbox{outside $f\inv(\cT\setminus\cT^{0})$.}
	\end{array}
 	\right.
	\end{equation*}
Otherwise, say $n_{2}=0$, the subscripts $0$ of the products are set to be vacuous. 
And $f\cup g=f$.
Note that if either $n_{1}$ or $n_{2}$ is $0$, $\varpi_{\sT}^{n_{1},n_{2}}$ becomes the identity map.

Taking the union over all $n_{1},n_{2}\ge0$ to $\varpi_{\cT}^{n_{1},n_{2}}$, we obtain a natural isomorphism
	$$
	\varpi_{\cT}\colon
	(\Map(\pp^{1},\cT)\times \Map(\pp^{1},\cT))_{0}
	\cong \Map(\pp^{1},\cT)\times_{S^{\bullet}\cc}(S^{\bullet}\cc\times S^{\bullet}\cc)_{0}.
	$$

We apply Drinfeld's factorization property to the product stack $\cT\times\cT$.
Note that the substack $(\cT\times\cT)^{0}:=\cT^{0}\times\cT^{0}$ has the Cartier divisor complement $(\cT\times(\cT\setminus\cT^{0}))\cup((\cT\setminus\cT^{0})\times\cT)$. 
Thus we obtain the factorization isomorphism $\varpi_{\cT\times\cT}$.
On the other hand there is an obvious identification 
	\begin{equation}\label{eq: obvious identification}
	\Map(\pp^{1},\cT\times\cT)=\Map(\pp^{1},\cT)\times\Map(\pp^{1},\cT).
	\end{equation}
Under this identification, the two factorization isomorphisms are equal:
	\begin{equation}\label{eq: factorization isomorphisms}
	\varpi_{\cT\times\cT}=\varpi_{\cT}\times\varpi_{\cT}.
	\end{equation}
This will be used in the proof of Theorem \ref{th: tensor product commutes with factorization}  later on.

We give a more general rule for the factorization isomorphisms (Proposition \ref{prop: F A}).
Since this formulation will not be used in the proof of Theorem \ref{th: tensor product commutes with factorization}, the readers can skip safely Proposition \ref{prop: F A}.
Let $\cT'$ be a stack as above with an open substack $\cT'{}^{0}\cong \Spec\cc$ such that $\cT'\setminus\cT'{}^{0}$ is  a Cartier divisor.
We suppose that there is a morphism $\sA\colon\cT\to \cT'$ such that $\sA(\cT^{0})$ is a substack of $\cT'{}^{0}$. 
Then we have the induced morphism between the spaces of based maps
	$$
	\sA_{*}\colon\Map(\pp^{1},\cT)\to \Map(\pp^{1},\cT'),\quad f\mapsto \sA f.
	$$

\begin{prop}\label{prop: F A}
The following composites are identical morphisms on $(\Map(\pp^{1},\cT)\times \Map(\pp^{1},\cT))_{0}$
	$$
	(\sA_{*}\times \Id_{S^{\bullet}\cc\times S^{\bullet}\cc})\varpi_{\cT}=\varpi_{\cT'}(\sA_{*}\times\sA_{*}).
	$$ 
\end{prop}
\proof
First we prove the composite $\varpi_{\cT'}(\sA_{*}\times\sA_{*})$ is well-defined.
In other words if $f,g\in\Map^{\ge1}(\pp^{1},\cT)$ have the disjoint supports (i.e.\ $f\inv(\cT\setminus\cT^{0})\cap g\inv(\cT\setminus\cT^{0})=\emptyset$), so do $\sA f,\sA g$.
But this is obvious.

Now to prove the proposition, we need to check that the pasted functions in both sides are identical, i.e.\ $\sA(f\cup g)=\sA f\cup \sA g$.
This is clear from the definition of pasting.
\qed\vskip.3cm

Now we identify the instanton space with a space of maps following \cite[\S3.1]{BFG}.
The above $(\pp^{1},\infty)$ is denoted by $(\pp^{1}_{h},\infty_{h})$ from now on.
We use another copy $(\pp^{1}_{v},\infty_{v})$.
For an algebraic group $G$ with a compact subgroup $K$ such that $K_{\cc}=G$, let $\cM^{G}_{\pp^{1}_{v}}$ be the moduli stack of (holomorphic) principal $G$-bundles on $\pp^{1}_{v}$ with trivialization at $\infty_{v}$.
For a based map $f\colon(\pp^{1}_{h},\infty_{h})\to (\cM^{G}_{\pp^{1}_{v}},[G\times\pp^{1}_{v}])$, the pull-back via $f\times \Id_{\pp^{1}_{v}}$ of the universal bundle over $\cM^{G}_{\pp^{1}_{v}}\times\pp^{1}_{v}$ defines a principal $G$-bundle on $\pp^{1}_{h}\times\pp^{1}_{v}$ with trivialization along $(\pp^{1}_{h}\times\infty_{v})\cup (\infty_{h}\times\pp^{1}_{v})$.
This principal $G$-bundle over $\pp^{1}_{h}\times\pp^{1}_{v}$ with a trivialization can be regarded as a principal $G$-bundle over $\pp^{2}$ with a trivialization using the isomorphism $\pp^{2}\setminus l_{\infty}\cong (\pp^{1}_{h}\times\pp^{1}_{v})\setminus ((\pp^{1}_{h}\times\infty_{v})\cup (\infty_{h}\times\pp^{1}_{v}))$.
By identifying the principle $G$-bundles over $\pp^{2}$ with trivialization along $l_{\infty}$ as framed $K$-instantons, we obtain an isomorphism $\Map(\pp_{h}^{1},\cM^{G}_{\pp^{1}_{v}})\cong \cM^{K}$, where $\cM^{K}$ denotes the framed $K$-instanton space.
The reverse morphism is also naturally defined.

According to \cite[\S3.3]{BFG}, there exists a moduli stack $\tcT$ (a thick Grassmannian) with a natural morphism $\tcT\to \cM^{G}_{\pp^{1}_{v}}$ which induces an isomorphism $\Map(\pp^{1}_{h},\tcT)\cong \Map(\pp^{1}_{h},\cM^{G}_{\pp^{1}_{v}})$ (\cite[Proposition 3.4]{BFG}).
According to \cite[\S2.2]{BFG}, there is a pro-algebraic group action on $\tcT$ whose quotient stack, say $\cT$, is equipped with the one-point open subset $\cT^{0}$ parametrizing the trivial objects.
Furthermore  $\cT^{0}$ and the substack of $\cM^{G}_{\pp^{1}_{v}}$ of trivial $G$-bundles are pulled back to an identical substack of $\tcT$, say $\tcT^{0}$.
We have a natural isomorphism $\Map(\pp^{1}_{h},\tcT)\to \Map(\pp^{1}_{h},\cT)$.
Using the morphism $\sE_{\cT}\colon \Map(\pp^{1}_{h},\cT)\to S^{\bullet}\cc$, we define the degree $n$ component $\Map^{n}(\pp^{1}_{h},\tcT)$ and hence $\Map^{n}(\pp^{1}_{h},\cM^{G}_{\pp^{1}_{v}})$ as before. 

\begin{prop}\label{prop: Map space instanton space} $($Cf.\ \cite[\S2.3 and Lemma 3.2]{BFG}$)$ 
If $K$ is simple, for any $n\ge0$ there are natural isomorphisms
	$$
	\Map^{n}(\pp^{1}_{h},\cT)\cong\Map^{n}(\pp^{1}_{h},\tcT)\cong \Map^{n}(\pp^{1}_{h},\cM^{G}_{\pp^{1}_{v}})\cong \cM^{K}_{n}.
	$$ 
\end{prop}

Let $G=\SL(N)$ and $K=\SU(N)$.
We fix an instanton number $n$.
We interpret the morphism $\sE_{\cT}^{n}$ in terms of ADHM data $(B_{1},B_{2},i,j)\in \bM_{\cc^{n},\cc^{N}}$.
We denote by $E^{n}\colon\cM^{K}_{n}=\mu\inv(0)^{\reg}/\GL(n)\to S^{n}\cc$ the morphism sending $[(B_{1},B_{2},i,j)]\mapsto$ the set of eigenvalues of $B_{1}$ counted with multiplicities.

\begin{prop}\label{prop: E and Ek}
Under the identification $\Map^{n}(\pp^{1}_{h},\cT)\cong \cM^{K}_{n}$, we have
$\sE_{\cT}^{n}=E^{n}$.
\end{prop}

\proof
Let $f\in \Map^{n}(\pp^{1}_{h},\cT)$ and $\cV'_{f}$ be the corresponding vector bundle over $\pp^{1}_{h}\times\pp^{1}_{v}$.
Let $\cM'$ be the substack in $\cM_{\pp_{v}^{1}}^{G}$ of nontrivial framed vector bundles.
Then both $\cM'$ and the Cartier divisor $\cT\setminus \cT^{0}$ are pulled back to $\tcT\setminus\tcT^{0}$.
By the construction of $(\cT,\cT^{0})$, under the identification $\Map(\pp^{1}_{h},\cT)\cong \Map(\pp^{1}_{h},\cM^{G}_{\pp^{1}_{v}})$ the Cartier divisor $\sE^{n}_{\cT}(f)=f\inv(\cT\setminus\cT^{0})$ in $\pp^{1}_{h}$ has the support of the points $z\in\cc$ such that the restriction $\cV'_{f}|_{\{z\}\times\pp^{1}_{v}}$ is a nontrivial framed vector bundle.

Let $\cV_{f}$ be the framed vector bundle on $\pp^{2}$ corresponding to $\cV'_{f}$.
Let $\pp^{1}_{z}:=\{z_{1}=z\}$ in $\pp^{2}$.
It is clear that $\cV_{f}|_{\pp^{1}_{z}}\cong \cV'_{f}|_{\{z\}\times\pp^{1}_{v}}$ under the obvious identification $\pp^{1}_{z}=\{z\}\times \pp^{1}_{v}$.
In order to prove the proposition, we need to write $\cV_{f}|_{\pp^{1}_{z}}$ in terms of ADHM data as follows:
If we denote by $(B_{1},B_{2},i,j)\in\mu\inv(0)^{\reg}$ an ADHM datum of $\cV_{f}$, Barth's correspondence (\cite[(2.6)]{Lecture}) identifies $\cV_{f}$ as the cohomology sheaf $\Ker(\beta)/\Image(\alpha)$ with a natural trivialization.
To be precise, let $\cO:=\cO_{\pp^{2}}$ for short.
Let $z_{1},z_{2}$ be the affine coordinates of $\cc_{h}\times\cc_{v}=\pp^{2}\setminus l_{\infty}$ where $\cc_{h}=\pp^{1}_{h}\setminus\infty_{h}$ and $\cc_{v}=\pp^{1}_{v}\setminus\infty_{v}$.
Let $0\to \cO(-1)^{\oplus n}\stackrel\alpha\to \cO^{\oplus n}\oplus\cO^{\oplus n}\oplus \cO^{\oplus N}\stackrel\beta\to \cO(1)^{\oplus n}\to 0$ be a monad, where $\alpha,\beta$ are matrix forms in $z_{1},z_{2}$ given by 
	$$
	\alpha=\begin{pmatrix} 
	B_{1}-z_{1}
	\\
	B_{2}-z_{2}
	\\
	j
	\end{pmatrix}
	,\quad
	\beta=\begin{pmatrix} 
	-B_{2}+z_{2}
	&
	B_{1}-z_{1}
	&
	i
	\end{pmatrix}.
	$$
The natural trivialization on $\Ker(\beta)/\Image(\alpha)$ is induced from the projection $\cO_{l_{\infty}}^{\oplus n}\oplus\cO_{l_{\infty}}^{\oplus n}\oplus \cO_{l_{\infty}}^{\oplus N}\to \cO_{l_{\infty}}^{\oplus N}$. 

Let $\alpha_{z},\beta_{z}$ be the morphisms induced from $\alpha,\beta$ by restriction to $\pp^{1}_{z}$.
The framed vector bundle $\cV_{f}|_{\pp^{1}_{z}}$ is given as $\Ker(\beta_{z})/\Image(\alpha_{z})$ (with a natural trivialization).
Now by the above discussion the proposition amounts to the equivalence: $\Ker(\beta_{z})/\Image(\alpha_{z})$ is a nontrivial framed vector bundle if and only if $B_{1}-z$ is non-invertible.

We take the restriction of the above monad to $\pp^{1}_{z}$.
The framed vector bundle on $\pp^{1}_{z}$ induced by restriction is trivial if and only if there is an isomorphism $\cO_{\pp^{1}_{z}}^{\oplus N}\cong \Ker(\beta_{z})/\Image(\alpha_{z})$.
Note that any morphism $\cO_{\pp^{1}_{z}}^{\oplus N}\to \Ker(\beta_{z})/\Image(\alpha_{z})$ is induced from a morphism $\cO_{\pp^{1}_{z}}^{\oplus N}\to \Ker(\beta_{z})$.
This fact follows from taking the ext functor $\Ext^{\bullet}(\cO_{\pp^{1}_{z}}^{\oplus N},-)$ to the short exact sequence $0\to \cO_{\pp^{1}_{z}}(-1)^{\oplus n}\stackrel{\alpha_{z}}\to \Ker(\beta_{z})\to \Ker(\beta_{z})/\Image(\alpha_{z})\to0$ and then using $\Ext^{1}(\cO_{\pp^{1}_{z}}^{\oplus N},\cO_{\pp^{1}_{z}}(-1)^{\oplus n})=0$ (the Serre duality).
This argument also shows that any isomorphism $\cO_{\pp^{1}_{z}}^{\oplus N}\cong \Ker(\beta_{z})/\Image(\alpha_{z})$ which restricts to the natural trivialization on the infinity point $\infty_{z}$ is unique.
Hence an isomorphism $\cO_{\pp^{1}_{z}}^{\oplus N}\cong \Ker(\beta_{z})/\Image(\alpha_{z})$, if any, is uniquely given by $e_{m}\mapsto (u_{m},v_{m},e_{m})$ for some $u_{m},v_{m}\in\cc^{n}$  ($m=1,2,...,N$), where $e_{m}$ is the $m^{\mathrm{th}}$ standard basis element of $\cc^{N}$ and $(u_{m},v_{m},e_{m})$ is a solution of the equation $\beta_{z}(u_{m},v_{m},e_{m})=0$, i.e.\ $(-B_{2}+z_{2})u_{m}+(B_{1}-z)v_{m}+i(e_{m})=0$.
This equation assures that $u_{m}=0$ by the degree reason in $z_{2}$ and moreover that $v_{m}$ uniquely exists if and only if $B_{1}-z$ is invertible.
Therefore $\Ker(\beta_{z})/\Image(\alpha_{z})$ is a trivial framed vector bundle if and only if $B_{1}-z$ is invertible.

So far the proof of the proposition at set-theoretic level has been done, i.e.\ $\sE^{n}(f)=E^{n}(\cV_{f})$ as sets.
To prove the general case, it suffices to check that any Cartier divisor $\sE^{n}(f)$ is always a limit of Cartier divisors of the form $\sE^{n}(f)$ consisting of points with only multiplicities $1$.
This fact follows from that the image of $(\mu_{1}\inv(0))^{k}\times_{S^{(1^{k})}\cc}(S^{1}\cc)^{k}_{0}$ via the ADHM data version factorization morphism $\pi_{(1^{k})}$ is Zariski open dense in $\mu\inv(0)$ (see the proof of Proposition \ref{prop: k2}). 
\qed\vskip.3cm

Recall that as far as $\cM^{K}_{n}$ has the ADHM description, the map $E^{n}$ is defined using the ADHM data.
Thus the above proof also works for the classical groups:

\begin{cor}\label{cor: classical group En}
For any classical group $G$ and $n\ge0$, we have $\sE_{\cT}^{n}=E^{n}$ under the identification $\Map^{n}(\pp^{1}_{h},\cT)\cong \cM^{K}_{n}$.
\end{cor}

Therefore this corollary combined with Proposition \ref{prop: E and Ek} yields the following corollary on the two kinds of factorization isomorphisms:

\begin{cor}\label{cor: old and new factorization}
For any classical group $G$ and $n\ge0$, under the identification $\Map^{n}(\pp^{1}_{h},\cT)\cong \cM^{K}_{n}$, the  factorization isomorphisms $\varpi_{\cT}^{n_{1},n_{2}}$ and $\opi_{\eta}$ coincide, where $n=n_{1}+n_{2}$ and $\eta$ is the partition of length $2$ consisting of $n_{1},n_{2}$.
\end{cor}

Here when $G=\SO(4)$, $\cM^{K}_{n}$ denotes the disjoint union $\coprod_{n=n_{1}+n_{2}}\cM^{\SO(4,\rr)}_{(n_{1},n_{2})}$.

Now we are ready to prove Theorem \ref{th: tensor product commutes with factorization}.
\vskip.3cm

\textit{Proof of Theorem \ref{th: tensor product commutes with factorization}.} 
Let $\cT$ be the stack corresponding to $\cM_{\pp^{1}_{v}}^{\Sp(1)}$.
Recall the natural isomorphisms $\Map(\pp^{1}_{h},\cT\times\cT)\cong \cM^{\Spin(4)}\cong \cM^{\SO(4,\rr)}$.
Under the isomorphism the obvious identification \eqref{eq: obvious identification}: $\Map(\pp^{1}_{h},\cT)\times \Map(\pp^{1}_{h},\cT)=\Map(\pp^{1}_{h},\cT\times\cT)$ restricts to the tensor product isomorphism $T_{(n_{1},n_{2})}\colon\cM^{\USp(1)}_{n_{1}}\times \cM^{\USp(1)}_{n_{2}}\to \cM^{\SO(4,\rr)}_{(n_{1},n_{2})}$ in the diagram \eqref{eq: factorization diagram}.
Recall here that the tensor product morphism is originally defined over the moduli spaces of vector bundles corresponding to the instanton spaces (\S\ref{subsec: Geometry of the instanton moduli space}).
Recall also that the associated vector bundle to a $K$-instanton is defined as $P\times_{K}W$, where $W$ is the vector representation of $K$.
In our case, the tensor product of two $\USp(1)$-instantons $P_{1},P_{2}$ (or equally associated vector bundles) is nothing but the vector bundle $(P_{1}\times_{S^{4}}P_{2})\times_{\SO(4,\rr)}(W_{1}\otimes W_{2})$, where  $W_{1},W_{2}$ are (identical) vector representations of $\USp(1)$, because $W_{1}\otimes W_{2}$ becomes the vector representation of $\SO(4,\rr)$.

Now the commutativity of \eqref{eq: factorization diagram} is equivalent to \eqref{eq: factorization isomorphisms}, because the (horizontal) factorization rational maps in \eqref{eq: factorization diagram} can be rewritten in terms of $\varpi_{\cT}$ and $\varpi_{\cT\times\cT}$ due to Corollary \ref{cor: old and new factorization}.
This proves Theorem \ref{th: tensor product commutes with factorization}.
\qed\vskip.3cm


\section{Nilpotent symmetric matrices}\label{App: nilpotent symmetric matrices}

Let $V$ be a vector space of dimension $k$ with a symplectic or orthogonal form $(\,,\,)$.
Let $B$ be a nilpotent endomorphism in $\fp(V)$.
An aim of this appendix section is to construct some nice basis of $V$ as below.
Note that a similar basis for $\fg(V)$ is obtained by Springer and Steinberg \cite[\S{IV}.2.19]{SS}.

\begin{lem}\label{lem: G}
Let $B\in \fp(V)$ be a nilpotent endomorphism.

$(1)$
If $V$ is symplectic, there exist nonzero vectors $v_{l},v_{l}'$ $(1\le l\le s)$ and $d_{l}\in\zz_{\ge0}$ of $V$ such that $B^{a}v_{l},B^{a'}v_{l}'$ $(0\le a,a' \le d_{l})$ form a basis of $V$ and all the pairing among them are zero except
	$$
	(B^{a}v_{l},B^{d_{l}-a}v_{l}')=1.
	$$

$(2)$
If $V$ is orthogonal, there exist nonzero vectors $u_{l}$ $(1\le l\le t)$ and $e_{l}\in\zz_{\ge0}$ of $V$ such that $B^{b}u_{l}$ $(0\le b\le e_{l})$ form a basis of $V$ and all the pairing among them are zero except
	$$
	(B^{b}u_{l},B^{e_{l}-b}u_{l})=1.
	$$
\end{lem}

\proof
We sketch the proof.
Let us define a bilinear form $|\,,\,|$ on $\Image B$ as $|Bv,Bv'|:=(v,Bv')=(Bv,v')$.
Then it is a symplectic (resp.\ orthogonal) form on $\Image B$ if so is $(\,,\,)$.
It is clear that $B|_{\Image B}$ is symmetric with respect to $|\,,\,|$.

If $B=0$, the statement itself amounts to the existence of symplectic or orthogonal basis.
Otherwise, by the induction hypothesis we can take $v_{l},v_{l}'$ or $u_{l}$ such that their $B$-images satisfy the statement for $(\Image B,|\,,\,|)$.
We set $s=t=1$ as the general case comes from orthogonal decomposition.
Then $d_1,e_1\ge1$.
Now direct computation shows $v_{1},v_{1}'-(v_{1},v_{1}')B^{d_{1}}v_{1}'$ or $u_{1}-\frac12(u_{1},u_{1})B^{e_{1}}u_{1}$ satisfy the lemma.
\qed\vskip.3cm

\begin{cor}
\label{cor: G''}
Let $B\in \fp(V)$.
Then
$G(V).B=\GL(V).B\cap \fp(V)$.
\end{cor}

\proof
It is clear that $G(V).B\subset \GL(V).B\cap \fp(V)$.
We prove the opposite inclusion.
Let $A\in \GL(V).B\cap \fp(V)$.
Let $a_1,a_2...,a_e$ be the mutually distinct eigenvalues of $B$.
They are also eigenvalues of $A$.
Let $V_l(A),V_l(B)$ be the generalized $a_l$-eigenspaces of $A,B$ respectively where $l=1,2,...,e$.
Then the isomorphism induced by the bijection between the bases in Lemma \ref{lem: G} for nilpotent endomorphisms $A|_{V_l(A)}-a_l\Id_{V_l(A)}$ and $B|_{V_l(B)}-a_l\Id_{V_l(B)}$ gives conjugacy equivalence of $A,B$.
\qed\vskip.3cm

A partition $\eta:=(\eta_1,\eta_2,...)$ is of \textit{even type} if  $\heta_{l}$ is even for all $l$ where $\heta=(\heta_{1},\heta_{2},...)$ denotes the dual partition of $\eta$.
Let $B$ be a nilpotent endomorphism in $\gl(V)$.
Then $\dim \Ker B^{l}/\Ker B^{l-1}$, $l\ge1$, form a partition for $l\ge1$.
Its dual partition is called the \textit{associated partition} of $B$.
Even if $B$ is not nilpotent we also define the associated partitions as follows:
Let $a_1,a_2,...,a_e$ be the mutually distinct eigenvalues of $B$ and $V_l$ be the generalized $a_l$-eigenspace of $B$.
By the \textit{associated partitions} of $B$, we mean the associated partitions of the nilpotent endomorphisms $B|_{V_l}-a_l\Id_{V_l}$ ($l=1,2,..,e$).

Lemma \ref{lem: G} (1) yields the following corollary:

\begin{cor}
\label{cor: G}
Let $B\in \fp(V)$.
If $V$ is symplectic, all the associated partitions are of even type.
\qed
\end{cor}

The following lemma gives the converse of Corollary \ref{cor: G}.

\begin{lem}
\label{lem: H}
Let $e\ge1$ and $\eta^1,\eta^2,...,\eta^e$ be partitions with $|\eta^1|+|\eta^2|+\cdots+|\eta^e|=k$.

$(1)$
If $V$ is orthogonal, there exists $B\in\fp(V)$ having $\eta^1,\eta^2,...,\eta^e$ as the associated partitions.

$(2)$
If $V$ is symplectic and moreover each partition $\eta^{l}$ is of even type, there exists $B$ in $\fp(V)$ having $\eta^1,\eta^2,...,\eta^e$ as the associated partitions.
\qed
\end{lem}


\end{document}